%% file: ncube1fin.tex
\def\version{May 6, 2003}

\documentclass[12pt]{article}
\usepackage[psamsfonts]{amsfonts}
\usepackage{amsmath,amssymb,amsthm}
\usepackage[dvips]{graphicx}
\usepackage{eepic}
\input def99c 
\UseSection  
\renewcommand{\to}{\rightarrow}

\newcommand{\Pro}{{\mathbb P}_p}

\newcommand{\Var}{{{\rm Var}_p}}
\newcommand{\Exp}{{{\mathbb E}_p}}
\newcommand{\Expn}{{{\mathbb E}_{p_n}}}

\newcommand{\ben}{\begin{enumerate}}
\newcommand{\een}{\end{enumerate}}

\newcommand{\prob}[2][]{{{\mathbb P}_{#1}(#2)}}

\newcommand{\Cmax} {{\Ccal}_{\rm max}}

\newcommand{\shift}   {\!\!\!\!}

\newcounter{countC}  
\newcounter{countR}  
\newcommand{\gr}{{\mathbb G}}
\newcommand{\ver}{{\mathbb V}}
\newcommand{\edg}{{\mathbb B}}
\newcommand{\qn}{{\mathbb Q}_n}
\newcommand{\egr}\Exp
\newcommand{\cn}{\Omega}

\newcommand{\constt}{a_0}
\newcommand{\constp}{a_1}
\newcommand{\constopup}{a_2}
\newcommand{\constoplow}{a_3}
\newcommand{\epO}{{\epsilon_0}}

\newcommand{\Z}{\Zbold}

\newcommand{\conn}{\leftrightarrow}

\newcommand{\dbc}{\Leftrightarrow}
\newcommand{\ct}[1]  { \stackrel{#1}{\conn} }

\newcommand{\ctx}[1]  {\leftarrow\shift\!\xrightarrow{#1}}
\newcommand{\nc}  { \conn  {\hspace{-2.5ex} /} \hspace{1.8ex}   }

\newcommand{\AND}       {\;\&\;}

\newcommand{\bigo}{O}

\newcommand{\Expg}{{\mathbb E}_{p,\gamma}}
\newcommand{\Prog}{{\mathbb P}_{p,\gamma}}
\newcommand{\Expcg}{{\mathbb E}_{p_c,\gamma}}

\title  {Random subgraphs of finite graphs: \\
        I. The scaling window under the triangle condition
        }

\author{Christian Borgs%
\thanks{Microsoft Research, One Microsoft Way, Redmond,
WA 98052, USA. {\tt borgs@microsoft.com},
{\tt jchayes@microsoft.com}}
\and
Jennifer T.\ Chayes$^*$
\and
Remco van der Hofstad%
\thanks{Department of Mathematics and Computer Science,
Eindhoven University of Technology, P.O.\ Box  513,
5600 MB Eindhoven, The Netherlands.
{\tt rhofstad@win.tue.nl}}
\and
Gordon Slade%
\thanks{Department of Mathematics, University of British Columbia,
Vancouver, BC V6T 1Z2, Canada. {\tt slade@math.ubc.ca}}
\and
Joel Spencer\thanks{Department of Computer Science,
Courant Institute of Mathematical Sciences,
New York University, 251 Mercer St., New York, NY 10012, U.S.A.
{\tt spencer@cs.nyu.edu}
}}

\date\version

\begin{document}

\maketitle


\begin{abstract}
We study random subgraphs of an arbitrary finite
connected transitive graph
$\mathbb G$ obtained by independently deleting edges with probability
$1-p$.  Let $V$ be the number of vertices in $\gr$,
and let $\cn$ be their degree. We define the critical threshold
$p_c=p_c(\mathbb G,\lambda)$ to be the value of $p$ for
which the expected cluster size of a fixed vertex
attains the value $\lambda V^{1/3}$, where $\lambda$ is fixed and
positive. We show that for any such model, there is a phase
transition at $p_c$ analogous to the phase transition for the
random graph, provided
that a quantity called the triangle diagram
is sufficiently small at the threshold $p_c$.
In particular, we show that the largest cluster inside a scaling
window of size $|p-p_c|=\Theta(\cn^{-1}V^{-1/3})$ is of size
$\Theta(V^{2/3})$, while below this scaling window, it is much
smaller, of order $O(\epsilon^{-2}\log(V\epsilon^3))$, with
$\epsilon=\cn(p_c-p)$.
We also obtain an upper bound
$O(\cn(p-p_c)V)$ for the expected size of the largest cluster
above the window.
In addition, we define and analyze the
percolation probability
above the window and show that it is of order
$\Theta(\cn(p-p_c))$.  Among the models for which the triangle diagram
is small enough to allow us to draw these conclusions are the
random graph, the $n$-cube and certain Hamming cubes, as well as the
spread-out $n$-dimensional torus for $n>6$.
\end{abstract}



\section{Introduction and results}
\label{sec-intro}

\subsection{Background}
\label{sec-background}

Random subgraphs of finite graphs are of central interest
in modern graph theory.  The best known example is the
random graph $G(V,p)$.  It is defined as the subgraph
of the complete graph on $V$ vertices obtained by
deleting edges independently with probability $1-p$,
and was first studied by
Erd\H{o}s and R\'enyi in 1960 \cite{ER60}.  They showed
that when $p$ is scaled as
$(1+\epsilon)V^{-1}$, there is a phase transition
at $\epsilon = 0$ in the sense that the size of the largest
component is $\Theta (\log V)$ for $\epsilon < 0$, $\Theta (V)$
for $\epsilon > 0$, and has the nontrivial
behavior $\Theta (V^{2/3})$ for $\epsilon = 0$.

The results of Erd\H{o}s and R\'enyi were substantially
strengthened by Bollob\'as \cite{Boll84} and {\L}uczak
\cite{Lucz90}. In particular, they showed that the model has a
scaling window of width $V^{-1/3}$ in the sense that if $p =
(1 + \Lambda_V V^{-1/3})V^{-1}$, then the size of the largest
component is $\Theta (V^{2/3})$ whenever $\Lambda_V$ remains
uniformly bounded in $V$, is
$o(V^{2/3})$ whenever $\Lambda_V \rightarrow -\infty$,
and $\omega(V^{2/3})$
whenever
$\Lambda_V \rightarrow \infty$.

Considerably less is known for random subgraphs of other finite graphs.
An interesting example is the $n$-cube $\qn$, which has
vertex set $\{0,1\}^n$ and an edge joining any two vertices that
differ in exactly one component.  Let $V=2^n$ denote the number
of vertices in $\qn$.  It is
known since the work of Ajtai, Koml\'os and Szemer\'edi
\cite{AKS82} that for $p$ of the form $p=(1+\epsilon)n^{-1}$, the
largest component is of size $\bigo(n)$ when $\epsilon$ is fixed
and negative, and is of size at least $c2^n$ for some
positive $c = c(\epsilon)$ if $\epsilon$ is fixed and positive.
However, very little is known about the scaling window.
The best results
available are those of Bollob\'as, Kohayakawa and \L uczak in
\cite{BKL92}, who showed the following.
We use the standard terminology that a sequence of events $E_n$ occurs
{\em asymptotically almost surely} (a.a.s.)\ if
$\lim_{n\to\infty}\Pbold(E_n)=0$.
In \cite{BKL92}, it is shown that
that for $p=(n-1)^{-1}(1+\epsilon)$ the
size of the largest cluster is at most $O(n\epsilon^{-2})$ if
$\epsilon< -e^{-o(n)}$, is a.a.s. $(2\log 2)n\epsilon^{-2}(1+o(1))$
if $\epsilon \leq -(\log n)^2 (\log \log n)^{-1} n^{-1/2}$,
and is a.a.s. $2\epsilon 2^n(1+o(1))$ if $\epsilon
\geq 60 n^{-1}(\log n)^3$. Note that the resulting bounds,
while much sharper than those established in \cite{AKS82},  are
still far from establishing the behavior one would expect by
analogy with the random graph, namely a window of width
$\Theta(V^{-1/3})$ where the largest cluster is of size
$\Theta(V^{2/3})$, with different behavior outside the window on either side.

For random subgraphs of
finite subsets of
$\Z^n$, Borgs, Chayes, Kesten and Spencer
\cite{BCKS01} systematically developed a relationship
between critical exponents and the width of the scaling window.  In
particular, they determined
the size of the largest
component
inside, below, and above a suitably defined
window, under certain scaling and hyperscaling
hypotheses
(proved in $n = 2$ and conjectured to be
valid whenever $n\leq 6$).
These results gave
the appropriate version of the Erd\H{o}s and R\'enyi \cite{ER60},
Bollob\'as \cite{Boll84} and {\L}uczak
\cite{Lucz90} results for random subsets of $\Z^2$.


Very recently, there have been attempts to extend the Erd\H{o}s
and R\'enyi \cite{ER60} analysis to more general finite graphs.
Frieze, Krivelevich and Martin \cite{FKM02} showed that for random
subgraphs of pseudorandom graphs of $V$ vertices, there is a phase
transition in which the largest component goes from $\Theta(\log
V)$ to $\Theta(V)$.  Alon, Benjamini and Stacey \cite{ABS02} use
the methods of \cite{AKS82} to study the critical value for the
emergence of the $\Theta(V)$ component in random subgraphs of
general finite graphs of large girth. Note, however, that in the
language of the discussion above, both \cite{FKM02} and
\cite{ABS02} consider only $\epsilon$ {\em fixed}; they do not get
any results on the scaling window.

In this paper, we study conditions under which
random subgraphs of arbitrary finite graphs behave like the random
graph $G(V,p)$,
both with respect to the critical point and the
scaling window.  More precisely, let $\gr$ be a finite
connected  {\em
transitive} graph with $V$ vertices of degree $\cn$.  Consider
random subgraphs of $\gr$ in which edges are deleted independently
with probability $1-p$.  We show that if
a finite version of
the so-called
triangle diagram is sufficiently small at a suitably
defined transition point $p_c$ (see below), then the model behaves like
the random graph in the sense that inside a window of width
$|p-p_c|=\Theta(\cn^{-1}V^{-1/3})$ the largest cluster is of order
$\Theta(V^{2/3})$ while it is of order $o(V^{2/3})$ below this
window.
These results are essentially optimal within and below the scaling
window.  While we do obtain results much stronger than previous
results above the scaling window, our bounds in this region are still
far from optimal.  It is likely that a condition beyond the triangle
condition (e.g., an expansion condition) will be necessary to achieve
optimal results in this region above the window.

For percolation on {\em infinite} graphs, the triangle diagram has
been recognized as an important quantity since the work of Aizenman
and Newman \cite{AN84}  who identified the so-called {\em triangle
condition} as a sufficient condition for mean-field behavior for
percolation on $\Z^n$.  Here, the term {\em mean-field
behavior} refers to the critical behavior of percolation on a
tree, which is well understood.  The triangle condition is defined
in terms of the {\em triangle diagram}
    \eq
    \lbeq{tria-def}
    \nabla_p(x,y)
    = \sum_{u,v \in \ver}
    \tau_p(x,u) \tau_p(u,v) \tau_p(v,y),
    \en
where the sum goes over the vertices of the underlying graph, and
$\tau_p(x,y)$ denotes the probability that $x$ and $y$ are
joined by a path of occupied edges (in the random subgraph
language, $\tau_p(x,y)$ is the probability that $x$ and $y$ lie in
the same component of the random subgraph). On $\Z^n$, the
triangle condition is the statement that at the threshold $p_c$,
 $\nabla_{p_c}(x,x)$ is finite.  The triangle
condition was proved on $\Z^n$ by Hara and Slade \cite{HS90a,HS94}, using the
lace expansion, for the nearest-neighbor model for $n \geq 19$
and for a wide class of spread-out (long-range) models for
 all $n>6$.

Let $\chi(p)$ denote the expected size of the cluster containing a fixed
vertex.
Aizenman and Newman used a differential inequality for
$\chi(p)$ to show that the triangle
condition implies that as $p\nearrow p_c$, the expected cluster
size diverges like $(p_c-p)^{-\gamma}$ with $\gamma=1$.
Subsequently, Barsky and Aizenman \cite{BA91} showed, in particular, that the
triangle condition also implies that as $p\searrow p_c$ the
percolation probability goes to zero like $(p-p_c)^{{\beta}}$ with
$\beta=1$.  Their proof is based on differential inequalities for the
{\em magnetization}.  These inequalities,
which were motivated by an earlier inequality of
Chayes and Chayes \cite{CC86d,CC87b}, had been used previously by
Aizenman and Barsky to prove
sharpness of the percolation phase transition on $\Z^n$ \cite{AB87}.
The exponents $\gamma$ and ${\beta}$ are examples of
{\em critical exponents}.  For percolation on a tree, the above
behavior for the percolation probability and the expected cluster
size can be relatively easily established with $\gamma = {\beta}
=1$.

In order to apply the above methods to prove mean-field behavior for
percolation on
finite graphs, several hurdles must be overcome. The first is
the fact that it is {\it a priori} unclear how even to define the
critical value $p_c$.  Second, the triangle
condition must be modified, since $\nabla_p(x,y)$ is always finite on
a finite
graph.  Third, the method of integration of the differential inequalities
of \cite{AB87,AN84,BA91} requires
that at $p_c$, the expected cluster size diverges, which is again
not possible on a finite graph $\gr$. All these
facts, which we deal with below,
require substantial modification and generalizations of the
methods and concepts of \cite{AB87,AN84,BA91}.

In addition to the methods involving differential inequalities,
our results are based on a second
set of techniques, developed in \cite{BCKS01}, relating critical exponents
and the width of the scaling window.
We will apply these methods here
to obtain information on the size of the largest cluster from
information on the cluster-size distribution.

The results of this paper are valid assuming the triangle condition.
For the complete
graph $G(V,p)$, we will easily verify the triangle condition below,
thereby reproducing some of the known results for the phase transition
in the random graph.
In \cite{BCHSS04b}, we will use the lace expansion to verify the
triangle condition for several other examples of finite graphs,
including the $n$-cube and various tori with vertex set $\{0,1,\ldots,
r-1\}^n$.  This leads to several new results for these models;
see Section~\ref{sec-examples} below.

\subsection{The setting}
\label{sec-setting}

Let $\mathbb G=(\mathbb V,\mathbb B)$ be a finite graph.
The vertex set $\mathbb V$ is any finite set, and
the set of bonds (or edges) $\mathbb B$ is a subset of the set
of all two-element subsets $\{x,y\}\subset\mathbb V$.
The {\em degree} of a vertex $x\in\mathbb V$ is
the number of bonds containing $x$.   A bijective map
$\varphi:\mathbb V\to \mathbb V$ is called a {\em graph isomorphism} if
$\{\varphi(x),\varphi(y)\}\in\mathbb B$ whenever
$\{x,y\}\in\mathbb B$, and $\mathbb G$ is called {\em transitive} if
for every pair of vertices $x,y\in \mathbb V$ there is a
graph-isomorphism $\varphi$ with $\varphi(x)=y$.  Transitive graphs
are by definition regular, i.e., each vertex has the same degree.

Let $\gr$ be an arbitrary finite, connected, transitive graph
with $V$ vertices of degree $\cn$.  We study percolation on $\gr$,
in which each of the bonds is {\it occupied} with probability $p$
independently of the other bonds, and {\it vacant} otherwise. We
denote probabilities and expectations in the resulting product
measure by $\Pro(\cdot)$ and $\Exp(\cdot)$, respectively.

As usual, we say that $x$ {\it is connected to} $y$, written as $x
\conn y$, when there is a path from
$x$ to $y$ consisting of occupied bonds.  We define the {\it
connectivity function} $\tau_p(x,y)$ by
\eq
\lbeq{tau-def}
\tau_p(x,y)=\Pro(x\leftrightarrow y).
\en
We denote by $C(x)$ the cluster of a vertex $x$, that is, the
set of all vertices in $\gr$ which are connected to $x$,
and by $|C(x)|$ the number of vertices in this cluster.
Note that the distribution of $|C(x)|$ is invariant
under the automorphisms of $\gr$, and hence independent
of $x$.  Instead of $|C(x)|$, we will therefore
often study $|C(0)|$, where $0$, the ``origin'', is
an arbitrary fixed vertex in $\ver$.

Our main results involve the
{\em cluster size distribution},
    \eq\lbeq{Pgeqk}
    P_{\geq k}(p)
    =\Pro(|C(0)| \geq k),
    \en
the {\it susceptibility}
\eq\lbeq{chidef}
    \chi(p) =\Exp|C(0)|,
    \en
(i.e., the expected size of the cluster of a fixed vertex), and
the maximal cluster size
  \eq\lbeq{Cmax}
    |\Cmax| = \max\{|C(x)| :  x \in \gr  \}.
    \en
By definition, the function
$\chi$ is strictly monotone increasing on the interval $[0,1]$,
with $\chi(0)=1$ and $\chi(1)=V$.  Also,
\eq
\lbeq{chitau}
    \chi(p) = \Exp \sum_{x \in \mathbb V} I[x \in C(0)]
    = \sum_{x \in \mathbb V} \tau_p(0,x).
\en

Recall that for $G(V,p)$
the largest cluster inside the transition window is of order
$V^{2/3}$.  It is not difficult---in fact, easier---to
determine the expected cluster size inside the
window, which turns out to be of order $V^{1/3}$.
Motivated by this fact, we define the {\em critical threshold}
$p_c=p_c(\gr,\lambda)$ of a finite graph $\gr$ to
be the unique solution to the equation
    \eq
    \lbeq{pcdef}
    \chi(p_c) =  \lambda V^{1/3},\,
    \en
with $\lambda >0$.
There is some flexibility in the choice of $\lambda$, connected with
the fact that the transition takes place within a window and not at
a particular value of $p$.  A convenient choice is to take
$\lambda$ to be constant (independent of $V$). We will always assume that
$1<\lambda V^{1/3}<V$ so that  $p_c$ is well defined and
$0<p_c<1$.

The definition \refeq{pcdef} is appropriate for graphs that
obey mean-field behavior, which we expect only for graphs that are
in some sense ``high-dimensional.''  As we discuss in more detail
in Section~\ref{sec-uh}, a different definition of the critical threshold
would be appropriate for a graph providing a finite approximation
to $\Z^n$ for $n<6$.


\subsection{Main results}
\label{sec-mainresults}

In this section, we state our main results, which hold for
arbitrary finite connected transitive graphs, provided the triangle
diagram \refeq{tria-def} at $p_c$ is sufficiently small.
To be more precise, we will assume that
 \eq
\lbeq{tria-con}
    \nabla_{p_c}(x,y) \leq\delta_{x,y}+\constt
\en
for a sufficiently small constant $\constt$, a condition
we call the {\em finite-graph triangle condition}, or more briefly,
the {\em triangle condition}.
Although we have not done the necessary computations,
the constant $\constt$ need not
be extremely small, and we expect
our results to hold for $\constt$ of the order
of $\frac {1}{10}$.

By \refeq{chitau}, $\sum_{y \in \mathbb V} \nabla_{p}(x,y) = \chi(p)^3$.
As a consequence, the triangle condition implies
that
    \eq
    \lbeq{lambda-bd}
    \lambda^3\leq \constt + V^{-1}.
    \en
In other words, small $\lambda$ is a necessary condition for the
triangle condition to hold.  It turns out that it also
sufficient for many graphs $\gr$. For the random graph, this is shown in
Section~\ref{sec-RG}, and for several other models
in \cite{BCHSS04b}; see Section~\ref{sec-examples}.
Indeed, we will show that for
these models,
    \eq \lbeq{better-tria-con} \nabla_{p}(x,y)
    =\delta_{x,y}+ O(\cn^{-1}) + O(\chi^3(p)/V)
    \en
whenever
$\chi^3(p)/V$ is small enough.

Our results concerning the critical threshold are given in the
following theorem.  In its statement, we make the abbreviations
    \eq\lbeq{epsilon0}
    \epO = \frac 1{\chi(p_c)}=\lambda^{-1} V^{-1/3},
    \en
and
\eq
\lbeq{nabbardef}
    \bar\nabla_p = \max_{\{x,y\}\in \mathbb B}\nabla_{p}(x,y).
\en

\begin{theorem}[Critical threshold]
\label{main-thm-cv}
For all finite connected transitive graphs $\gr$, the following statements hold.

\noindent
i) If $\lambda >0$ and the triangle condition \refeq{tria-con} holds
for some $\constt<1$, then
    \eq\lbeq{omegapc-1}
    1-\epO
    \leq \cn p_c \leq
    \frac{1-\epO}{1-\constt}.
    \en

\noindent ii)
Given $0<\lambda_1 <\lambda_2 < \infty$, let $p_i$ be defined by
$\chi(p_i)=\lambda_i V^{1/3}$ ($i=1,2$).  If $\bar\nabla_{p_2}<1$, then
\eq
\lbeq{p2p1}
    \frac{\lambda_2 -\lambda_1}{\lambda_1\lambda_2 } \frac{1}{V^{1/3}}
    \leq
    \cn(p_2-p_1 )
    \leq
    \frac{1}{1-\bar\nabla_{p_2}}
    \frac{\lambda_2 -\lambda_1}
    {\lambda_1\lambda_2 } \frac{1}{V^{1/3}}.
\en
\end{theorem}

For example, if $\gr$ is the complete graph on $n$ vertices
(so that $V=n$ and $\cn = n-1$) and $p_2$ is
inside the transition window, then $p_1$ remains within the
transition window for any constant $\lambda_1<\lambda_2$.

Our
results concerning the subcritical phase are given in the
following theorem.

\begin{theorem}[Subcritical phase]
\label{main-thm-sub} There is a (small) constant $b_0>0$ such that
the following statements hold for all positive $\lambda$, all
finite connected transitive graphs $\gr$ and all $p$ of the form
$p=p_c-\cn^{-1}\epsilon$ with $\epsilon\geq 0$.

\noindent
i) If the triangle condition \refeq{tria-con} holds
for some $\constt<1$, then
    \eq\lbeq{chibd}
    \frac{1}{\epO+\epsilon}
    \leq
    \chi(p)
    \leq \frac{1}%
    {\epO+[1-\constt] \epsilon  }.
    \en

\noindent
ii) If the triangle condition holds for some
$\constt\leq b_0$ and if $\lambda V^{1/3}\geq b_0^{-1}$,
then
    \eq\lbeq{cmaxbd1}
    10^{-4} \chi^2(p)
    \leq
    \Exp\Big(|\Cmax|\Big)
    \leq
    2\chi^2(p)\log(V/\chi^3(p)),
    \en
\eq\lbeq{cmaxbd2}
    \Pro\Big( 
    |\Cmax|\leq 2\chi^2(p)\log(V/\chi^3(p))\Big)
    \geq 1-
    \frac {\sqrt{e}}{[2\log(V/\chi^3(p))]^{3/2}},
\en
and, for $\omega \geq 1$,
\eq
\lbeq{cmaxbdom}
    \Pro\Big(|\Cmax|\geq \frac{\chi^2(p)}{3600\omega}\Big)
    \geq\big(1+\frac {36\chi^3(p)}{\omega V}\Big)^{-1}.
\en
\end{theorem}

Our next theorem states our results inside the scaling
window.

\begin{theorem}[Critical Window]
\label{main-thm-critical}

Let $\lambda>0$ and $\Lambda<\infty$.
Then there are finite
positive constants $b_1,\dots,b_8$ such that the
following statements hold for all finite connected transitive graphs
$\gr$ provided the triangle-condition \refeq{tria-con} holds for
some constant $\constt\leq b_0$ and $\lambda V^{1/3}\geq b_0^{-1}$,
with $b_0$ as in Theorem~\ref{main-thm-sub}.
Let $p=p_c +\cn^{-1}\epsilon$
with $|\epsilon|\leq\Lambda V^{-1/3}$.

\noindent
i) If
$k\leq b_1 V^{2/3}$, then
\eq\lbeq{clszdis}
\frac{b_2}{\sqrt k}
\leq
P_{\geq k}(p)
\leq
\frac{b_3}{\sqrt k}.
\en

\noindent
ii) 
\eq\lbeq{LCEBd1-win}
    {b_4}V^{2/3}
    \leq
    \Exp\big[|\Cmax|\big]
    \leq
    {b_5}V^{2/3}
    \en
and, if $\omega\geq 1$, then
    \eq\lbeq{LCBd1-win}
    \Pro\Big(
    \omega^{-1} V^{2/3}\leq |\Cmax|\leq \omega V^{2/3}
        \Big)
    \geq 1-\frac{b_6}\omega.
    \en

\noindent
iii)
    \eq\lbeq{chiasy-win}
    b_7 V^{1/3}
    \leq \chi(p)\leq
    b_8 V^{1/3}.
    \en

\noindent In the above statements, the constants $b_2$ and $b_3$
can be chosen independent of $\lambda$ and $\Lambda$, the constants
$b_5$ and $b_8$ depend on $\Lambda$ and not $\lambda$, and the constants
$b_1$, $b_4$, $b_6$ and $b_7$ depend on both
$\lambda$ and $\Lambda$.

\end{theorem}

Our results on the supercritical phase are given in the following
theorem.

    \begin{theorem} [Supercritical phase]
    \label{main-thm-sup}

Let $\lambda>0$.  The following statements hold for all finite
connected transitive graphs $\gr$ provided the triangle-condition
\refeq{tria-con} holds for some constant $\constt\leq b_0$ and
$\lambda V^{1/3}\geq b_0^{-1}$, with $b_0$ as in Theorem~\ref{main-thm-sub}.
Let $p=p_c+\epsilon \cn^{-1}$ with $\epsilon \geq 0$.
\smallskip

\noindent i)
        \eq\lbeq{bound1-cmax}
        \Exp(|\Cmax|)
        \leq  21\epsilon V+ 7V^{2/3},
        \en
and, for all $\omega >0$,
\eq\lbeq{cmax.2A}
    \Pro\Big(|\Cmax|\geq \omega (V^{2/3}+\epsilon V)\Big)
\leq \frac{21}\omega .
    \en

\noindent ii)
    \eq\lbeq{chiasysup}
    \chi(p)
    \leq 81(V^{1/3}+\epsilon^2V).
    \en
\end{theorem}

Note that Theorem~\ref{main-thm-sup} does not give lower bounds
on the size of the largest supercritical cluster.
We believe that this is not a mere technicality.
Indeed, the formation of a giant component in the
random graph is closely related to the fact that
moderately large clusters have a significant chance
to merge into a single, giant component as
$\epsilon$ is increased beyond the critical value
by an amount of order $V^{-1/3}$.
This fact involves the geometry of the random
graph, and may not be true for arbitrary
transitive graphs obeying the triangle condition.
It would be interesting to know whether there
exists a sequence of transitive graphs
$\gr_n$ such that the largest cluster above the
window is $o(\epsilon V)$, at least if
$\epsilon V^{-1/3}\to\infty$ sufficiently slowly.
On the other hand, as we explain in more detail in
Section~\ref{sec-examples} below, our results apply to
the $n$-cube $\qn$, and for $\qn$ we prove complementary
lower bounds to the upper bounds of Theorem~\ref{main-thm-sup}
in \cite{BCHSS04c}.  Our proof of these upper bounds is valid for
$\epsilon \geq e^{-cn^{1/3}}$, and not in the full domain
$\epsilon \gg V^{-1/3}=2^{-n/3}$ where we would conjecture that
they are valid. The methods of \cite{BCHSS04c} rely heavily on
the specific geometry of $\qn$ and do not apply at the level of
generality of Theorem~\ref{main-thm-sup}.

We close this section with a theorem that
gives a more precise bound on the susceptibility below
the window, under the assumption that the stronger
triangle condition \refeq{better-tria-con} holds.
We make the constants in \refeq{better-tria-con} explicit by
assuming that
    \eq \lbeq{stronger-tria-con}
    \nabla_{p}(x,y)
    \leq\delta_{x,y}+ K_1\cn^{-1} + K_2\frac{\chi^3(p)}V
    \en
for some constants $K_1$, $K_2<\infty$ and all $p\leq p_c$.
Let
\eqalign
\lbeq{adef}
    a&=K_1\cn^{-1}+K_2\lambda^3,    \\
\lbeq{K2tildef}
    \tilde K_2  &=K_2/(1-a) ,\\
\lbeq{atildef}
    \tilde a(\epsilon)
    & =K_1\cn^{-1}+K_2\lambda^3\frac{\epO}{\epO+(1-a)\epsilon}.
\enalign

\begin{theorem}[Sharpened bounds]
\label{stronger-mainthm}
Let $\lambda>0$ and let \/ $\gr$ be a finite connected transitive graph
such that
\refeq{stronger-tria-con}
holds for all $p\leq p_c$, with the constant in
\refeq{adef} obeying $a<1$.  Let
$p=p_c-\cn^{-1}\epsilon$ with $\epsilon \geq 0$, and let
$\tilde{K}_2$ and $\tilde{a}(\epsilon)$ be given by
\refeq{K2tildef}--\refeq{atildef}.  Then
\eqalign
\lbeq{omegapc-1-strong}
{1-\epO}
&\leq
{\cn p_c}
\leq
\frac {1-\epO}{1-K_1 \cn^{-1}-\tilde K_2\lambda^3 \epO}, \\
    \lbeq{chibd-strong}
    \frac{1}{\epO+\epsilon}
    &\leq
    \chi(p)
    \leq \frac 1{\epO+[1-{\tilde a}(\epsilon)] \epsilon}.
    \enalign
\end{theorem}

The inequality \refeq{omegapc-1-strong} implies that $|\cn p_c-1|
=O(\cn^{-1}+\lambda^{-1}V^{-1/3})$.  The significance of \refeq{chibd-strong}
is most apparent if we consider a sequence of graphs with
$\lambda >0$ fixed, $V\to \infty$
and $\cn \to \infty$, for $\epsilon$ such that $\epsilon/\epO \to \infty$.
In this limit, \refeq{chibd-strong} implies that
\eq
    \chi(p) = \frac{1}{\epsilon}[1+o(1)].
\en
We will apply Theorem~\ref{stronger-mainthm} to the random graph in
Section~\ref{sec-RG}.

\subsection{General sequences of finite graphs}

To illustrate our theorems, it is instructive to consider
a sequence of finite connected transitive graphs $\gr_n=(\ver_n,\edg_n)$
with $|\ver_n|\to\infty$.  We will say that such a sequence
obeys the {\em finite-graph triangle condition} if there
exist a $\lambda>0$ such that the condition \refeq{tria-con}
holds for all $n$, with a constant $\constt$ that is at most as
large as the constant $b_0$ in Theorem \ref{main-thm-sub}.

Consider thus a sequence of finite connected transitive graphs
$\gr_n$ satisfying the finite-graph triangle condition.
Consider also a sequence of probabilities of the form
    \eq
    \lbeq{form-of-pn}
    p_n=p_c+\Lambda_n\cn^{-1}V_n^{-1/3}.
    \en
Motivated by the random graph (and our theorems)
we say that the sequence $p_n$ is {\em inside the window},
if $\limsup_{n\to\infty}|\Lambda_n|<\infty$, {\em below the window}
if $\Lambda_n\to -\infty$, and {\em above the window} if
$\Lambda_n\to \infty$ as $n\to\infty$.  In order to avoid
dealing with higher order corrections in
$\epsilon_n=\Lambda_n V_n^{-1/3}$, we assume here that
$\epsilon_n\to 0$.

Consider first a sequence {\em below the window}, i.e., assume that
$\Lambda_n\to -\infty$ as $n\to\infty$.  The first statement of
Theorem~\ref{main-thm-sub} then implies that
    \eq
    \chi(p_n)=\Theta(\epsilon_n^{-\gamma})
    \en
with $\gamma=1$, while
the second implies that
    \eq
    \Theta(\Lambda_n^{-2} V_n^{2/3})
    \leq
    \Expn\Big(|\Cmax|\Big)
    \leq
    \Theta(\Lambda_n^{-2} \, V_n^{2/3}\log\Lambda_n),
    \en
and
    \eq
    \phantom{\quad\text{a.a.s. as $n\to\infty$}}
    \Theta(\Lambda_n^{-2} V_n^{2/3})
    \leq
    |\Cmax|
    \leq
    \Theta(\Lambda_n^{-2}\, V^{2/3} \log\Lambda_n)
    \quad\text{a.a.s. as $n\to\infty$.}
    \en
Note that this implies,
in particular, that below the window,
$|\Cmax|=o(V_n^{2/3})$ a.a.s. as $n\to\infty$.

Next, consider a sequence $p_n$
{\em inside the window}, i.e., a sequence of the form
\refeq{form-of-pn} with $\limsup|\Lambda_n|<\infty$.
Theorem~\ref{main-thm-critical} then implies that
\eq
\chi(p_n)=\Theta(V_n^{1/3}),
\en
    \eq
    \Exp\big[|\Cmax|\big]=\Theta(V_n^{2/3}),
    \en
with the probability of the event
    \eq
    \omega(n)^{-1}
    \leq
    \frac{ |\Cmax|}{\Exp\big[|\Cmax|\big]}
    \leq
    \omega(n)
    \en
going to one whenever $\omega(n)\to\infty$ as $n\to\infty$.

Let us finally consider a sequence $p_n$
{\em above the window}, i.e., a sequence of the form
\refeq{form-of-pn} with $\Lambda_n\to\infty$.
Theorem~\ref{main-thm-sup}~i)~then implies that the
expected size of the largest cluster is
$\bigo(\epsilon_n V_n)$, and Theorem~\ref{main-thm-sup}~ii)~shows
that $\chi(p_n) = O(\epsilon_n^2 V_n)$.

\subsection{The percolation probability and magnetization}

It is a major result for percolation on $\Z^n$ that the value of $p$
for which $\chi(p)$ becomes infinite is the same as the value of
$p$ where the percolation probability, or order parameter, ${\mathbb
P_p}(|C(0)|=\infty)$, becomes positive
\cite{AB87,Mens86}. In the
present setting, since the graph is finite, there can be no infinite
cluster and the definition of the order parameter needs to be
adapted.  A natural
definition of the {\it finite-size order parameter}
is the ratio of the expected maximal cluster
size to the volume V:
   \eq\lbeq{theta}
     \theta(p)=
     \frac{\Exp(|\Cmax|)}V.
   \en

However, we are unable to prove a good lower bound on \refeq{theta}
in the supercritical regime,
and we therefore consider an alternative
definition in terms of the cluster size distribution
$P_{\geq k}(p)$.
Parameterizing $p$ as $p=p_c+\epsilon \cn^{-1}$,
we define the percolation probability by
    \eq\lbeq{theta_a}
    \theta_\alpha(p) =
    \Pro(|C(x)| \geq N_\alpha)
    = P_{\geq N_\alpha}(p),
    \en
where
    \eq\lbeq{N_a}
    N_\alpha=
    \frac 1{\epsilon^2}
    \big(\epsilon V^{1/3}\big)^\alpha.
    \en
Here
$\alpha$ is a constant with
$0<\alpha<1$.
This definition is motivated by the known behavior of the
random graph.  Above the window (corresponding to $\epsilon
V^{1/3}\to\infty$), it is known that
a.a.s., the largest component has size $|\Cmax| =
2 \epsilon  V[1+o(1)]$, while the second largest has size
$2\epsilon^{-2}\log(\epsilon^3 V)(1+o(1))$.  For the random graph
above the window, the cutoff $N_\alpha$ in \refeq{theta_a} is
therefore much larger than the second largest, and much smaller
than the largest cluster.  As a consequence, the ratio of
$\theta_\alpha(p)$ and $\theta(p)$ goes to one when considered on the
random graph above threshold.
(The above reasoning actually suggests the wider
range $0<\alpha<3$ for $\alpha$, but
for technical reasons
we require $0<\alpha<1$.)
Our results for $\theta_\alpha(p)$ are stated in the following theorem.

    \begin{theorem} [The percolation probability]
    \label{main-thm-pp}

Let $\lambda>0$ and $0<\alpha<1$.  Then there are finite
positive constants $b_9$, $b_{10}$, $b_{11}$, $b_{12}$
such that the following statements hold for all finite
connected transitive graphs $\gr$ provided the triangle-condition
\refeq{tria-con} holds for some constant $\constt\leq b_0$ and
$\lambda V^{1/3}\geq b_0^{-1}$, with $b_0$ as in Theorem~\ref{main-thm-sub}.
Let $p=p_c+\epsilon \cn^{-1}$.
\smallskip

\noindent i)
        \eq\lbeq{bound1}
         b_{10}\epsilon
         \leq
        \theta_\alpha (p)
        \leq
        27\epsilon,
        \en
where the lower bound holds when $b_9 V^{-1/3}\leq \epsilon\leq 1$
and the upper bound holds when $\epsilon V^{1/3} \geq 1$.

\noindent ii)
If $\max\{b_{12}V^{-1/3}, V^{-\eta}\}\leq
\epsilon\leq 1$, where
$\eta=\frac 13\frac{3-2\alpha}{5-2\alpha}$,
then
    \eq
    \lbeq{LCbdsup}
    \Pro\Big(
    |\Cmax|\leq [1 + (\epsilon V^\eta)^{-1}]\theta_\alpha(p) V
    \Big)
    \geq 1- \frac{b_{11}}{(\epsilon V^\eta)^{3-2\alpha}}.
    \en

\noindent In the above statements, the constants  $b_9$,
$b_{10}$, $b_{11}$ and $b_{12}$ depend on both $\alpha$ and
$\lambda$.
\end{theorem}

Theorem~\ref{main-thm-pp}~i) is analogous to results proved for percolation
on $\Z^n$ (assuming high $n$ for the upper bound) in \cite{AB87,BA91,HS90a}.
Theorem~\ref{main-thm-pp}~ii) shows that it is unlikely that the largest
supercritical cluster is larger than $\theta_\alpha(p)V$, at least for
$\epsilon$ not too small.
As we will describe in more detail in Section~\ref{sec-examples} below,
it is shown in \cite{BCHSS04c} that when $\gr$ is the $n$-cube,
it is possible also to prove a lower bound on $|\Cmax|$, so that
$|\Cmax|$ is of the same order of magnitude as $\theta_\alpha(p)V$,
at least when $\epsilon$ is not too small.
The fact that $\theta_\alpha(p)$ can be used in this way  serves
as further justification for the definition \refeq{theta_a}.
In \cite{BCKS01}, a similar approach was used for $\Z^n$ in low
dimensions.

Our analysis of $\theta_\alpha(p)$, and more generally our analysis
of the cluster
size distribution $P_{\geq k}(p)$, is primarly based on
an analysis of the magnetization.
Let $P_k(p)$ be the probability that the $|C(0)|=k$.
The
{\em magnetization} $M(p,\gamma)$ is defined by
    \eq
    M(p,\gamma)
    = 1 - \sum_{k=1}^{V} (1-\gamma)^{k} P_k(p).
    \en
Thus $M(p,\gamma)$ is essentially
the generating function for the sequence $P_k(p)$, and $M(p,0)=0$ for all $p$.
Estimates on $M(p,\gamma)$ for small $\gamma$ can be converted into
estimates on $P_k(p)$ for large $k$, via an analysis reminiscent of
a Tauberian theorem.  The name ``magnetization'' is used because $M(p,\gamma)$
is analogous the the magnetization in spin systems, and
the variable $h \geq 0$ defined by $\gamma = 1 -e^{-h}$ plays
the role of an external magnetic field in
that context.  Our main results for the magnetization are summarized
in the following theorem.

    \begin{theorem} [The magnetization]
    \label{main-thm-mag}
Assume that $\constt$ is sufficiently small, and let $0 \leq \gamma \leq 1$.

\noindent
i) If $p\leq p_c$ then
        \eq
        \lbeq{Mgeqthm}
        \frac 13
    \min\{\sqrt\gamma,\gamma\chi(p)\} \leq
    M(p,\gamma)
    \leq
     \min\{\sqrt{12\gamma},\,\gamma\chi(p)\}   .
    \en

\noindent
ii)
If $p=p_c+\cn^{-1}\epsilon$ and $\epsilon \geq 0$ then
    \eq
    \lbeq{Mleq1thm}
    M(p,\gamma)\leq
    \sqrt{12\gamma}+13\epsilon.
    \en
Let  $0\leq \alpha<1$ and $\rho>0$.
There is a positive $c=c(\alpha, \lambda)$ and
$b_{13} = b_{13}(\alpha ,\lambda, \rho)$ such that if
$b_{13}V^{-1/3} \leq \epsilon \leq 1$ then
    \eq
    \lbeq{Mgeq7thm}
    M(p,\rho N_\alpha^{-1})
    \geq
        c\epsilon  \min\{1,\rho^{1/(2-\alpha)}\} .
     \en
\end{theorem}

\subsection{Guide to the paper}

In Section~\ref{sec-ex}, we discuss several examples where our
general results can be applied.
In Section~\ref{sec-overview}, we indicate some of the main ideas
that enter into the proofs of our main results.


The following table indicates where the various theorems are proved.
The notation
[u.b.] refers to the upper bounds on $|\Cmax|$ and [l.b.] to the lower bounds.
\begin{center}
\begin{tabular} { | c || c | c | c | c | c | c |}   \hline
    Theorem
    & \ref{main-thm-cv}
    & \ref{main-thm-sub} i), ii) [u.b]
    & \ref{main-thm-sub} ii) [l.b]
    & \ref{main-thm-critical} i)
    & \ref{main-thm-critical} ii-iii)
    & \ref{main-thm-sup} i)
    \\ \hline
    Section \rule {0mm}{4mm}
        & \ref{sec-subp}
        & \ref{sec-subp}
        & \ref{sec-lbsc}
        & \ref{sec-proof}
        & \ref{sec-overview}
        & \ref{sec-proof}
        \\  \hline \hline
    Theorem
    & \ref{main-thm-sup} ii)
    & \ref{stronger-mainthm}
    & \ref{main-thm-pp} i)
    & \ref{main-thm-pp} ii)
    & \ref{main-thm-mag}
    &
        \\ \hline
    Section \rule {0mm}{4mm}
        & \ref{sec-chi-sup}
        & \ref{sec-subp}
        & \ref{sec-proof}
        & \ref{sec-supper-upb}
        & \ref{sec-M}
        &
        \\ \hline
\end{tabular}
\end{center}

There is no dependence on Section~\ref{sec-subp} in
Sections~\ref{sec-M}--\ref{sec-supper-upb}.
The bounds on the magnetization proved in Section~\ref{sec-M}
are crucial for Sections~\ref{sec-proof}--\ref{sec-supper-upb}.
Section~\ref{sec-lbsc} depends on Section~\ref{sec-proof},
which in turn depends on Section~\ref{sec-M}.  Sections~\ref{sec-chi-sup}
and \ref{sec-supper-upb} each depend on Section~\ref{sec-M} and on no
other section.  Sections~\ref{sec-lbsc}, \ref{sec-chi-sup} and
\ref{sec-supper-upb} are mutually independent.
Three differential inequalities, needed in Sections~\ref{sec-subp},
\ref{sec-M}, and \ref{sec-chi-sup},
are proved in Appendix~\ref{sec-app}.

\section{Examples}
\label{sec-ex}

\subsection{The random graph}
\label{sec-RG}

In this section, we illustrate both the finite-graph
triangle condition and our
results when $\gr$ is the random graph on $n$ vertices. In the
notation of the last section, we thus consider the graph
$\gr=K_n$, the complete graph on $n$ vertices, with $V=n$
vertices of degree $\cn=n-1$.

\subsubsection{The triangle condition for the random graph}
\label{sec-RGtri}

For the random graph, the triangle diagram can be explicitly
and easily calculated in terms of the expected cluster size $\chi(p)$,
as follows.
Due to the high degree of symmetry of the complete
graph, the two-point function takes on only the two distinct values
$\tau_p(x,x)=1$ and $\tau_p(x,y) = \tau$ (say) for $x\neq y$,
so that $\tau_p(x,y)=\delta_{x,y} + \tau(1-\delta_{x,y})$.
The triangle diagram \refeq{tria-def} is therefore given by
\eq
 \nabla_{p}(x,y)=
\begin{cases}
1+3(n-1)\tau^2+(n-1)(n-2)\tau^3
&\text{ if } x=y,
\\
3\tau+ 3(n-2)\tau^2 +
[1+(n-1)(n-2)]\tau^3
&\text{ if } x\neq y.
\end{cases}
\en
Also, by \refeq{chitau},
$\tau=(n-1)^{-1}(\chi(p)-1)$.  Since $\chi(p) \leq n$, this implies
that $\tau \leq n^{-1}\chi(p)$.  It is then straightforward to see
that
\eq\lbeq{stronger-tria-con-RG}
\nabla_{p}(x,y)
\leq
\delta_{x,y} +\frac{\chi^3(p)}n\left[1+3\chi^{-1}(p)+3\chi^{-2}(p)\right]
\leq \delta_{x,y}+7n^{-1}\chi^3(p).
\en
Recalling that by definition,
\eq
\lbeq{pc-def-RG}
\chi(p_c)=\lambda n^{1/3},
\en
we have thus obtained the triangle condition
\refeq{tria-con}
with $\constt= 7\lambda^3$.
In addition, \refeq{stronger-tria-con} holds with $K_1=0$ and $K_2=7$.

\subsubsection{The phase transition for the random graph}

Having verified the triangle condition, we can now apply
the results of Section~\ref{sec-mainresults} provided we take
$\lambda$ to be a sufficiently small constant.
Starting with Theorem~\ref{stronger-mainthm},
since $\cn = n-1=n(1+O(n^{-1}))$,
\refeq{omegapc-1-strong} implies that
\eq
p_c=\frac 1n(1+\bigo(n^{-1/3})).
\en
While we cannot expect that $p_c=1/n$ (in fact,
\refeq{omegapc-1-strong} implies that
$p_c<1/n$ if $\lambda$ is small enough), it
differs from the traditional value by only a small amount,
small enough to keep it inside the scaling window.
Thus our definition of $p_c$ is quite sensible for
the random graph.

In Theorems~\ref{main-thm-sub}--\ref{stronger-mainthm}, we have
used the parameter $\epsilon = \cn (p-p_c)$.  For the random graph,
we will use the scaling $p = p_c(1+\Lambda_n n^{-1/3})$, which
corresponds to $\epsilon = np_c \Lambda_n n^{-1/3}$.  Then up to constants,
$\epsilon$  is equivalent to
\eq
    \epsilon_n = \Lambda_n n^{-1/3}.
\en
Note that if $\Lambda_n \to -\infty$ then
for $K_1=0$ we have $\tilde{a}(\epsilon) = \Theta (|\Lambda_n^{-1}|)$, and
\refeq{chibd-strong} implies the simpler statement
\eq
    \chi(p) = \frac{1}{|\epsilon_n|}(1+O(\Lambda_n^{-1})),
\en
as claimed below in \refeq{chibd-RG}.

The conclusions of Theorems~\ref{main-thm-sub}--\ref{main-thm-pp}
for the random graph are summarized in the following theorem.

\begin{theorem}
\label{RG-thm}
Let $p=p_n=p_c(1+\Lambda_n n^{-1/3})$
with $p_c=p_c(n,\lambda)$ defined by \refeq{pc-def-RG}.
There exists a constant $\lambda_0$ such that the
following statements are true for all fixed,
strictly positive $\lambda\leq \lambda_0$, with
the constants implicit in our $\bigo(\cdot)$
and $\Theta(\cdot)$ possibly depending on
$\lambda$.

\noindent
i) (Subcritical phase). If $\Lambda_n\to-\infty$ as $n\to\infty$ then
    \eq\lbeq{chibd-RG}
    \chi(p)=\frac {n^{1/3}}{|\Lambda_n|}
    (1+O(\Lambda_n^{-1})),
    \en
    \eq\lbeq{cmaxbd1-RG}
    n^{2/3}\Theta(\Lambda_n^{-2})
    \leq
    \Exp\Big(|\Cmax|\Big)
    \leq
    6n^{2/3}\,
    \frac{\log|\Lambda_n|}{\Lambda_n^{2}}
    {(1+\bigo(\Lambda_n^{-1}))}
   ,
    \en
with
    \eq\lbeq{cmaxbd2-RG}
    n^{2/3}\Theta(\Lambda_n^{-2})
    \leq
    |\Cmax|
    \leq
    6n^{2/3}\,
    \frac{\log|\Lambda_n|}{\Lambda_n^{2}}
    {(1+\bigo(\Lambda_n^{-1}))}
    \quad \text{a.a.s.\ as $n \to \infty$.}
    \en

\noindent
ii) (Critical window).
If $\Lambda=\limsup|\Lambda_n|<\infty$ as $n\to\infty$ then
\eq
\chi(p)=\Theta(n^{1/3}),
        \quad \quad
    \Exp\big[|\Cmax|\big]=\Theta(n^{2/3}),
    \en
with
    \eq
    \omega(n)^{-1}
    \leq
    \frac{ |\Cmax|}{\Exp\big[|\Cmax|\big]}
    \leq
    \omega(n)
    \quad \text{a.a.s.\ as $n \to \infty$}
    \en
whenever $\omega(n)\to\infty$ as $n\to\infty$.
If $kn^{-2/3}$ is small enough (depending on $\Lambda$), then
\eq
P_{\geq k}(p)=\Theta(k^{-1/2}).
\en

\noindent
iii) (Supercritical phase).
Let $0<\alpha<1$.  If $\Lambda_n\to\infty$ as $n\to\infty$
and $\epsilon_n=\Lambda_n n^{-1/3}\to 0$ then
\eq
\chi(p)=O( n^{1/3}\Lambda_n^2),
\quad \quad
\Exp(|\Cmax|)=O(\epsilon_n n),
\en
and
\eq
\theta_\alpha(p)=\Theta(\epsilon_n).
\en
If $\Lambda_n\to\infty$ at least as fast
as $n^{\eta}$,
where $\eta=\frac 13\frac {3-2\alpha}{5-2\alpha}$, then
\eq
\Pro\Big(|\Cmax|\leq
\big(1+\frac{1}{\epsilon n^{\eta}}\big)\theta_\alpha(p)n\Big)
\geq 1-\bigo\Big(\frac{1}{(\epsilon n^{\eta})^{3-2\alpha}}\Big).
\en
\end{theorem}

It is interesting to compare Theorem~\ref{RG-thm} with
previously known
results for the phase transition in the random graph.
Since
$p_c= n^{-1} + O(n^{-4/3})$, if we change our parametrization
to $p = n^{-1} + \Lambda n^{-4/3}$ then we effectively change
$\Lambda$ by
a constant.  This affects the constants in the critical
window and has an asymptotically negligible effect in the
subcritical and supercritical phases.  This new parametrization
is the standard parametrization
used (with $\lambda$ instead of $\Lambda$) in much of the random graph
literature.
We will refer to the book of Janson,
\L uczak and Rucinski \cite{JLR00}, where references to
the original literature can be found.
Results in \cite{JLR00} are expressed in terms of the
variable $s$, where $\frac{n}{2}+s$ gives the number of occupied
edges, and our formulas can be compared to theirs by setting
$s=\frac{\Lambda}{2}n^{2/3}$.

In the subcritical phase, we show that the largest component has size
between $c_1n^{2/3}\Lambda^{-2}$ and $c_2n^{2/3}\Lambda^{-2}\log |\Lambda|$,
while \cite[Theorem~5.6]{JLR00} gives, in particular,
that the largest component is asymptotically of size
$6n^{2/3}\Lambda^{-2}\log|\Lambda|$. The constant 6 in the upper bounds in
\refeq{cmaxbd1-RG}--\refeq{cmaxbd2-RG} is therefore sharp.

In the critical window, we show that the largest component has size
$\Theta(n^{2/3})$, while \cite[Theorem~5.20]{JLR00}
gives, in particular,  that the largest component
has size $X_1' n^{2/3}$ where $X_1'$ is a random variable with a nontrivial
distribution over $(0,\infty)$.

In the supercritical phase, let
$\epsilon = \Lambda n^{-1/3}$, so that $p=\frac{1}{n}(1+\epsilon)$.
We show that the largest component has size $O(\epsilon n)$.
As mentioned below the statement of Theorem~\ref{main-thm-sup},
we have no lower bound on the largest subcritical cluster in our
general setting.
In \cite[Theorem~5.12]{JLR00}, the
largest component is shown asymptotically to have size $2\epsilon n$
(their $\overline{s}$
is asymptotic to $s$ when $\epsilon\rightarrow 0$).
Moreover, \cite[Theorem~5.7]{JLR00} yields
that the $r$th  largest component (for any fixed $r\geq 2$)
has size asymptotic to $6n^{2/3}\Lambda^{-2}\log\Lambda$.
We are unable to get any reasonable upper bounds on the size of the
second largest component.

Although our results are not state-of-the-art for the random graph,
it is nevertheless striking that they follow from a general
theory that makes no calculation specific to the random graph
apart from the simple verification of the finite-graph triangle condition
in Section~\ref{sec-RGtri}.
More importantly, our theorems apply
much
more generally, to models such as the
$n$-cube and finite tori in $\Z^n$ for $n>6$, where they imply
strong
new results.

\subsection{The $n$-cube and several tori}
\label{sec-examples}

In \cite{BCHSS04b}, we use the lace expansion to
prove quite generally that for finite graphs
that are tori the triangle condition for
percolation is implied by a certain triangle condition for simple
random walk on the graph.  As we show in \cite{BCHSS04b},
the latter is easily verified for the
following graphs with vertex set $\{0,1,\ldots, r-1\}^n$:
\begin{enumerate}
\item
The narrow torus:
an edge joins vertices that differ by $1$ in exactly one component,
with the periodic boundary condition that $0$ and $r-1$ differ by 1,
for $r\geq 2$ fixed and $n\to \infty$.  For $r=2$, this is
the $n$-cube.
\item
The Hamming torus:
an edge joins vertices that differ in exactly one component, again
with the periodic boundary condition,
for $r\geq 2$ fixed and $n\to \infty$.
\item
The wide torus in high dimensions:
the same edge set as in (i) but now $n$ is large and fixed and we
study the limit $r \to \infty$ to approximate $\Z^n$.
\item
The wide spread-out torus in dimensions $n > 6$:
an edge joins vertices $x=(x_1,\ldots,x_n)$ and $y=(y_1,\ldots,y_n)$
if $\max_{i=1,\ldots,n}|x_i-y_i| \leq L$ (with periodic boundary conditions)
with $n > 6$ fixed, $L$ large and fixed, in the limit $r \to \infty$ to
approximate range-$L$ percolation on $\Z^n$.
\end{enumerate}
Our conclusions thus apply to the percolation
phase transition for each of the above examples.
The above examples are all high-dimensional graphs.
We do not expect the triangle condition to hold for low-dimensional
graphs, and in particular do not expect the triangle condition to
hold for the wide tori in dimensions $n \leq 6$.  Nor do we expect
the conclusions of our
theorems to hold in low dimensions.

Combined with \cite{BCHSS04b},
our results show that the phase transition for
percolation on the $n$-cube $\qn$ shares several features with
the phase transition for the random graph.
In particular, it follows from the triangle condition for $\qn$
proved in \cite{BCHSS04b} and Theorem~\ref{stronger-mainthm} that
$p_c(\qn,\lambda)= n^{-1}+O(n^{-2})$, for any sufficiently small choice of
$\lambda$.  In \cite{HS03a},
this series is substantially extended.
In \cite{BCHSS04c}, we use the lower bound on the percolation probability
of \refeq{bound1} to prove a lower bound on the largest supercritical
cluster for the $n$-cube.  This leads to a substantial improvement
of some of the results of \cite{AKS82,BKL92}.

Our results for the wide tori in high dimensions show that in
a window of width $r^{-n/3}$ centered at $p_c=p_c(r,n)$,
the largest cluster has size $\Theta (r^{2n/3})$.
It is interesting to compare this with a previous result for $\mathbb Z^n$.
For $p= p_c(\mathbb Z^n)$,
consider the restriction of percolation configurations to a large
box of side $r$, under the {\em bulk}\/ boundary condition in which
the clusters in the box are defined to be the intersection of the box with
clusters
in the infinite lattice (and thus clusters in the box need not be
connected within the box).
How large is the largest
cluster in the box, as $r \to \infty$?  The combined results of
Aizenman \cite{Aize97} and Hara, van der Hofstad and Slade \cite{HHS03}
show that for spread-out models with $n>6$
the largest cluster has size of order $r^4$,
and there are order $r^{n-6}$ clusters of this size.
For the nearest-neighbor model in dimensions $n \gg 6$, the same results
follow from the combined results of \cite{Aize97} and Hara \cite{Hara00}.
The size $r^4$ for the largest critical cluster
size is different than the $r^{2n/3}$ that we prove
for $p=p_c(r,n)$ under the {\em periodic} boundary condition of the torus.
Aizenman \cite{Aize97} had raised
the question whether a change from bulk to periodic boundary conditions
would change the $r^4$ to
$r^{2n/3}$.  It would be interesting to attempt to extend our results,
to show that $p_c(\mathbb Z^n)$ lies inside the critical window centered
at $p_c(r,n)$ for large $r$, thereby providing
an affirmative answer to Aizenman's question.

\section{Overview of the proofs}
\label{sec-overview}

\subsection{Differential inequality for the susceptibility}

The results for the critical threshold and the subcritical
susceptibility, stated in Theorems~\ref{main-thm-cv}, \ref{main-thm-sub}~i)
and \ref{stronger-mainthm} are all derived from the differential
inequality
\eq
\lbeq{dis}
    [1-\bar\nabla_p] \cn \leq - \frac{d\chi^{-1}(p)}{dp} \leq \cn,
\en
with $\bar\nabla_p$ defined by \refeq{nabbardef}.
This differential inequality was proved by Aizenman and Newman
\cite{AN84} with infinite graphs in mind, but its
proof applies also to finite transitive graphs.  We recall
the proof of \refeq{dis} in Appendix~\ref{sec-chidif}.
The triangle condition is used to
bound the left side of \refeq{dis} from below.
In Section~\ref{sec-revsub}, we will show that integration of
\refeq{dis} leads directly to proofs of
Theorems~\ref{main-thm-cv}, \ref{main-thm-sub}~i)
and \ref{stronger-mainthm}.

\subsection{Differential inequalities for the magnetization}

Aizenman and Barsky \cite{AB87}
used differential inequalities for the magnetization
to prove sharpness of the phase transition for percolation on $\Z^n$.
In \cite{BA91},
they derived a complementary differential inequality, assuming the
triangle condition, which implied that on $\Z^n$ the
magnetization and  percolation probability behave asymptotically
as $M(p_c,\gamma) = \Theta (\sqrt{\gamma})$ and
$\Pbold_p(|C(0)|=\infty) = \Theta (p-p_c)$.  In Section~\ref{sec-M},
we recall the statement of the differential inequalities of \cite{AB87},
and in
Appendix~\ref{sec-M-diff-ine} we derive a variant of the complementary
differential inequality of \cite{BA91}.  In Section~\ref{sec-M},
we show how to integrate
the differential inequalities to obtain the bounds on the magnetization
stated in Theorem~\ref{main-thm-mag}.  In performing the integration,
care is required to deal with the finite size effects.

The bounds on the magnetization proved in Theorem~\ref{sec-M} lie at
the heart of our method.  They play
a crucial role in all of Sections~\ref{sec-proof}--\ref{sec-supper-upb}
and in the proofs of Theorems~\ref{main-thm-sub}~ii)--\ref{main-thm-pp}.

\subsection{The cluster size distribution}

The magnetization is a generating function
for the sequence $P_k(p)$, and its behavior for small $\gamma$ is closely
related to the behavior of $P_{\geq k}(p)$ for large $k$.
This is made precise in Section~\ref{sec-proof}, where
Theorem~\ref{main-thm-critical}~i) and related bounds on $P_{\geq k}(p)$
are obtained from the bounds on the magnetization proven in
Section~\ref{sec-M}.  The upper bounds on the magnetization easily
lead to upper bounds on the cluster size distribution for all $p\in [0,1]$.
The lower bounds are more difficult.  We will need matching upper
and lower bounds on $M(p,\gamma)$ to obtain good
lower bounds on $P_{\geq k}(p)$, and, in the supercritical phase,
our lower bounds on $M(p,\gamma)$
are in the restricted form given in \refeq{Mgeq7thm},
with $\gamma$ proportional
to $N_\alpha^{-1}$.
Our bounds on
$P_{\geq k}(p)$ then lead to a proof of the bounds on
$\theta_\alpha (p) = P_{\geq N_\alpha}(p)$ stated in
Theorem~\ref{main-thm-pp}~i).

\subsection{The scale of the largest cluster}
\label{sec-LC-Heuristic}

\subsubsection{The random variable $Z_{\geq k}$}

Given $k>0$, let
\eq\lbeq{Zdef}
Z_{\geq k}=\sum_{x\in\mathbb V} I[|C(x)|\geq k]
\en
denote the number
of vertices that lie in clusters of size $k$ or larger.
Then
\eq\lbeq{ZExp}
\Exp (Z_{\geq k}) =V P_{\geq k}(p).
\en
By definition, $|\Cmax|\geq k$ if and only
if $Z_{\geq k}\geq k$, and hence, by the Markov inequality,
\eq
\begin{aligned}
\lbeq{cmax.3}
\Pro\big(|\Cmax|\geq k\big)
&
\leq\frac{VP_{\geq k}(p)}k.
\end{aligned}
\en
and
\eq
\begin{aligned}
\lbeq{cmax.4}
\Exp\big(|\Cmax|\big)
\leq k+\Exp\big(Z_{\geq k}\big)
=k+VP_{\geq k}(p).
\end{aligned}
\en
In addition,
\eq
\lbeq{Cmaxk}
|\Cmax|=\max\{k: Z_{\geq k}\geq k\},
\en
and hence the random variables $\{Z_{\geq k}\}_{k \geq 1}$
provide a characterization
of $|\Cmax|$.

\subsubsection{A useful heuristic}
\label{sec-uh}

The identity \refeq{Cmaxk} suggests that if the distribution of
$Z_{\geq k}$ is sufficiently concentrated
about its mean, then it should be the case that
\eq
\Exp(|\Cmax|)= \Theta\big( \max\{k: \Exp(Z_{\geq k})\geq k\} \big).
\en
Define
$k_0=k_0(p)$ to be the solution of the equation
\eq
k_0 = \Exp(Z_{\geq k_0}) = VP_{\geq k_0}(p).
\en
Then we are led to expect that
\eq
\lbeq{k0}
\Exp(|\Cmax|)=\Theta(k_0(p)).
\en
Under certain conditions, this heuristic was made rigorous in
\cite{BCKS01} to analyze percolation on finite subsets of $\Z^n$, $n \leq 6$,
and it underlies our approach to obtaining bounds on
$|\Cmax|$ from bounds on the cluster size distribution $P_{\geq k}(p)$.
As a reality check, we note that
for the random graph
it is not difficult to verify that as $|\epsilon|\to 0$,
\eq
\lbeq{k0rg}
k_0(p)=
\begin{cases}
2\epsilon^{-2}\log (\epsilon^3 V)(1+o(1))
&\text{below the window,}
\\
\Theta(V^{2/3})
&\text{inside the window,}
\\
2\epsilon V(1+o(1))
&\text{above the window.}
\end{cases}
\en
To leading order, this is precisely the size of
the largest cluster of the random graph, confirming \refeq{k0}.
Since we are working in settings where random graph scaling should
apply, \refeq{k0rg} also serves as a guide for our more general
transitive graphs.

In particular,
as noted in \cite{BCKS01},
if at the critical threshold we have
\eq
    P_{\geq k}(p_c) = \Theta (k^{-1/\delta}),
\en
then $k_0 = \Theta(V^{\delta/(\delta+1)})$ and \refeq{k0} predicts
that $\Ebold_{p_c}(|\Cmax|)=\Theta(V^{\delta/(\delta+1)})$.
This provides a connection between
the critical exponent $\delta$ and the size of the largest cluster
at criticality.  If we assume that
$\chi(p_c)$ is well approximated by $\Ebold_{p_c}(|\Cmax|) \Pbold_{p_c}(0 \in
\Cmax) \approx V^{\delta/(\delta+1)} V^{-1+\delta/(\delta+1)}$,
it also suggests that the correct definition of
the critical threshold, in general, is that value of $p$ for
which $\chi(p) = V^{(\delta - 1)/(\delta +1)}$.  Again, a constant
factor $\lambda$ could be introduced on the right side without
significant effect.
For a critical branching
process, it is the case that $\delta =2$.  For percolation on $\Z^n$
with $n$ sufficiently large,
it was proved in \cite{HS00b}
that $\delta = 2$ in the sense that
$P_k(p_c) = ck^{-3/2}(1+k^{-a})$ for some $a,c >0$.
On the other hand, it is believed
that $\delta$ is strictly greater than $2$
below the upper critical dimension $n=6$.  Thus we expect that the results
of Section~\ref{sec-examples} do not extend to wide tori for $n<6$,
and that our definition of $p_c$ also requires modification in this case,
namely in \refeq{pcdef}, the exponent $1/3$ should be replaced by
${(\delta - 1)/(\delta +1)}$.

We have in mind a high-dimensional graph $\gr$ for which cycles
are of limited importance.
Since each vertex has $\cn$ neighbors, criticality corresponds to
$p\cn \approx 1$, or $p_c \approx \cn^{-1}$.
According to the above, the value $\delta =2$ gives the familiar value
$V^{2/3}$ for the largest critical cluster.
How near to $p_c$ can we expect this behavior to hold, i.e., how
wide is the critical window?  Let us consider
$p<p_c$, which is easier.  If $p=p_c-\cn^{-1}\epsilon$,
we expect that a birth process with survival rate $1-\epsilon$
gives a good approximation, so that
\eq
P_{\geq k}(p)
\approx
\frac {{\rm const}}{\sqrt k}e^{-k\epsilon^2 /2}.
\en
The exponential is unimportant as long as
$\epsilon V^{1/3}\leq O(1)$, leading to $P_{\geq k}(p) \approx k^{-1/2}$
and thus $|\Cmax| \approx V^{2/3}$.  This suggests that
the system behaves critically when $\epsilon = O(V^{-1/3})$.

\subsubsection{Our method of proof}

Our proofs of bounds on $|\Cmax|$ proceed as follows.  For $p \leq p_c$,
we obtain an upper bound on $|\Cmax|$ by applying the upper bound
    \eq
    \lbeq{C0ub}
    P_{\geq k}(p) \leq \sqrt{\frac ek}e^{-k/(2\chi^2)},
    \en
which is valid for $k\geq \chi^2(p)$.  The bound \refeq{C0ub}
is proved in
\cite[Proposition 5.1]{AN84} and \cite[(6.77)]{Grim99}
(the proofs apply directly to any finite transitive graph).
We use \refeq{C0ub} in conjunction with \refeq{cmax.3}--\refeq{cmax.4},
choosing $k$ in accordance with the subcritical case in \refeq{k0rg}.
The details are carried out in Section~\ref{sec-scp1}.
For a lower bound on $|\Cmax|$, we prove a variance estimate for
$Z_{\geq k}$ and use this in conjunction with the second moment method.
The details are carried out in Section~\ref{sec-lbsc}.

Inside the critical window, our bounds on $|\Cmax|$ follow directly
from monotonicity and the subcritical and supercritical bounds.
This is discussed in Section~\ref{sec-red}.

In the supercritical phase,
the bounds on $|\Cmax|$ of Theorem~\ref{main-thm-sup}~i) follow directly
from our upper bounds on $P_{\geq k}$, and are derived in
Section~\ref{sec-proof}.  To prove the upper bound on $|\Cmax|$ stated
in Theorem~\ref{main-thm-pp}~ii),
we prove another variance
estimate for $Z_{\geq k}$.
This estimate allows us to bound the
probability that $Z_{\geq N_\alpha}$ differs from
its expectation $V\theta_\alpha(p)$ by more than a small multiple
of $V\theta_\alpha(p)$.
The variance of $Z_{\geq N_\alpha}$
is ultimately estimated in terms of the magnetization, and
the details are carried out in Section~\ref{sec-supper-upb}.
The restriction $\epsilon \geq V^{-\eta}$ in \refeq{LCbdsup}
(with $\eta \in (\frac{1}{9},\frac {3}{15})$ for $\alpha \in (0,1)$)
means that
this upper bound on $|\Cmax|$ has not yet been proven for all $p$
above the window.

\subsection{The supercritical susceptibility}
\label{sec-sucov}

The magnetization has a useful and standard probabilistic
interpretation. We define i.i.d.\ vertex variables taking the value
``green'' and ``not green'' by declaring that each $x \in \mathbb V$ is
green with probability $\gamma\in [0,1]$.  The vertex variables are
independent of the bond variables. Let $\Gcal$ denote the random
set of green vertices.  Then, by definition,
    \eq
    M(p,\gamma)
    =\sum_{k=1}^{V} [1-(1-\gamma)^{k}] \mathbb P_{p}(|C(0)|=k)
    = \Pbold_{p,\gamma}(0 \conn \Gcal),
    \en
where $\{0 \conn \Gcal\}$ denotes the event that $0 \conn x$ for some
$x \in \Gcal$.  Let
\eq
\lbeq{chiMpg}
    \chi(p,\gamma) = (1-\gamma) \frac{\partial}{\partial \gamma}
    M(p,\gamma)
    =
    \sum_{k=0}^{V} k(1-\gamma)^k \Pbold_{p} (|C(0)|=k)
    =
\Expg\big(|C(0)|I(0\nc \Gcal)\big)
\en
and
\eq
\lbeq{chiperp-rewz}
\chi_\perp(p,\gamma)
=
    \sum_{k=0}^{V} k[1-(1-\gamma)^k] \Pbold_{p} (|C(0)|=k)
    =
\Expg\big(|C(0)|I(0\leftrightarrow \Gcal)\big).
\en
The proof of Theorem~\ref{main-thm-sup}~ii) is based on the decomposition
\eq\lbeq{chisum1}
\chi(p)
=\chi(p,\gamma)+\chi_\perp(p,\gamma),
\en
which is valid for all $\gamma \in [0,1]$.

It follows from
\refeq{Mleq1thm}--\refeq{Mgeq7thm} (with $\alpha=0$ in the latter)
that $M(p,\epsilon^2) = \Theta(\epsilon)$
above the window.  For the
random graph the largest cluster above the window
has size of order $\epsilon V$, so
the origin is in the largest cluster with probability of order $\epsilon$.
Thus the probability that the origin is connected to the green set $\Gcal$
and the probability that the origin is in the largest cluster should
both be $\Theta(\epsilon)$, when we choose $\gamma = \epsilon^2$.
Thus we regard the green set $\Gcal$ as playing
the role of a kind of ersatz giant cluster, when $\gamma = \epsilon^2$.
From this perspective,
$\chi(p,\epsilon^2)$ corresponds to the expected cluster size
omitting the giant cluster, whereas $\chi_\perp(p,\epsilon^2)$
corresponds to
the expected cluster size of a vertex that is in the largest cluster.
Thus we might expect to prove that for $p \geq p_c$,
$\chi(p,\epsilon^2)$ is bounded above
by $O(\epsilon^{-1})$ while $\chi_\perp(p,\epsilon^2)$
is bounded above by $O(\epsilon^2 V)$.
An upper bound on $\chi(p,\gamma)$ will follow easily from
our bounds on the magnetization.
To obtain a bound of the form $O(\epsilon^2 V)$
for $\chi_\perp(p,\epsilon^2)$, we will make use of the random
variable
\eq\lbeq{ZGcal-def}
Z_{\Gcal}=
\sum_{x\in\ver}I(x\leftrightarrow \Gcal),
\en
which counts the number of vertices in clusters containing at least
one green vertex.   This will require a differential inequality
for the expectation of $Z_\Gcal^2$, which is proved in Appendix~\ref{sec-Zsup}.

\subsection{Proof of Theorem~\ref{main-thm-critical}~ii-iii)}
\label{sec-red}

Finally, we show that
the bounds of Theorem~\ref{main-thm-critical}~ii-iii) for the critical
window follow
from the bounds of
Theorems~\ref{main-thm-sub} and \ref{main-thm-sup}
for the subcritical and supercritical phases.

\smallskip \noindent
{\em Proof of \refeq{chiasy-win}.}
By the monotonicity of $\chi(p)$ in $p$,
the lower bound follows from the lower bound of \refeq{chibd}
(with $p=p_c-\Lambda \cn^{-1} V^{-1/3}$)
and the upper bound follows from the upper bound of \refeq{chiasysup}
(with $p=p_c+\Lambda \cn^{-1} V^{-1/3}$).
\qed

\smallskip \noindent
{\em Proof of \refeq{LCEBd1-win}.}
The upper bound follows from monotonicity of $\Ebold_p[|\Cmax|]$ in $p$
and the upper bound \refeq{bound1-cmax}
(with $p=p_c+\Lambda \cn^{-1} V^{-1/3}$).  The lower bound follows from
the lower bounds of \refeq{cmaxbd1} and \refeq{chiasy-win}.
\qed

\smallskip \noindent
{\em Proof of \refeq{LCBd1-win}.}
It follows from the upper bound of \refeq{LCEBd1-win} and Markov's
inequality that
\eq\lbeq{Cmax-up-a}
\Pro\Big(|\Cmax|\geq \omega V^{2/3}\Big)
\leq
\frac{b_5}{\omega}
\en
for all $\omega>0$.
For the complementary bound,
we bound $\Pbold_p(|\Cmax| \geq \omega^{-1}V^{2/3})$ below by
its value at $p=p_c-\Lambda V^{-1/3}$ and apply \refeq{cmaxbdom}
in conjunction with \refeq{chiasy-win}.
\qed

\section{The subcritical phase}
\label{sec-subp}

In Section~\ref{sec-revsub}, we apply a differential inequality
for $\chi(p)$ due to Aizenman and Newman \cite{AN84} to show that the
triangle condition \refeq{tria-con} implies the bounds
\refeq{omegapc-1}, \refeq{p2p1} and \refeq{chibd},
and that the stronger triangle condition \refeq{stronger-tria-con}
implies the bounds \refeq{omegapc-1-strong} and \refeq{chibd-strong}.
In Section~\ref{sec-scp1}, we apply the bound
\refeq{C0ub} on the cluster size
distribution, also due to \cite{AN84}, to
prove the upper bounds of \refeq{cmaxbd1}--\refeq{cmaxbd2}.

\subsection{The subcritical susceptibility and critical threshold}
\label{sec-revsub}

Recall from \refeq{nabbardef}
that $\bar\nabla_p = \max_{\{ x,y\} \in\edg}\nabla_p(x,y)$.
In Appendix~\ref{sec-chidif}, we prove the differential inequality
    \eq
    \lbeq{chi'ulz}
    [1-\bar\nabla_p] \cn \leq
    -\frac {d \chi^{-1}}{dp} \leq \cn,
    \en
which is valid for all $p \in (0,1)$.
The differential inequality
and its proof are due to Aizenman and Newman
\cite{AN84}.  Integration of \refeq{chi'ulz} over the interval $[p_1,p_2]$,
together with monotonicity of $\bar\nabla_p$ in $p$, gives
\eq
\lbeq{chiintegral}
    [1-\bar\nabla_{p_2}] \cn (p_2-p_1) \leq
    \chi^{-1}(p_1)-\chi^{-1}(p_2) \leq \cn (p_2-p_1).
\en

\smallskip \noindent
{\em Proof of \refeq{omegapc-1} assuming \refeq{tria-con}.}
We set $p_1=0$ and $p_2=p_c$ in \refeq{chiintegral}
and note that $\chi(0)=1$ and $\chi(p_c)^{-1} = \epsilon_0$,
to obtain \refeq{omegapc-1}.
\qed

\smallskip \noindent
{\em Proof of \refeq{p2p1}.}
This follows from \refeq{chiintegral} with
$p_i$ defined by $\chi(p_i) = \lambda_iV^{1/3}$.
\qed

\smallskip \noindent
{\em Proof of \refeq{chibd}.}
This follows from \refeq{chiintegral} with $p_1=p$ and $p_2=p_c$.
\qed

\smallskip \noindent
{\em Proof of \refeq{chibd-strong}.}
The lower bound has been proved already in \refeq{chibd}.
For the upper bound,
we first observe that the stronger triangle condition
\refeq{stronger-tria-con} implies \refeq{tria-con}
with $\constt=K_1\cn^{-1}+K_2\lambda^3$.  For $p\leq p_c$,
we may therefore use the upper bound of \refeq{chibd}
to see that
    \eq
    \lbeq{nabup}
    \nabla_{p}(x,y)
    \leq K_1\cn^{-1} + K_2\frac 1V\Big(\frac 1{\epO+(1-\constt)\epsilon}\Big)^3
    \en
for $x \neq y$.  We now integrate the lower bound of \refeq{chi'ulz} over
the interval $[p,p_c]$, using \refeq{nabup} to bound the triangle diagram.
This gives
\eqalign
\chi^{-1}(p)
-\epO
&\geq
\int_0^\epsilon d\tilde\epsilon
\Big[1-K_1\cn^{-1} -
K_2\frac 1V\Big(\frac 1{\epO+(1-\constt)\tilde\epsilon}\Big)^3
\big]
\nonumber
\\
&=\epsilon(1-K_1\cn^{-1})-\frac{K_2}{2V(1-\constt)}
\Big[\frac 1{\epO^2}
-\Big(\frac 1{\epO+(1-\constt)\epsilon}\Big)^2\Big]
\nonumber
\\
&=\epsilon(1-K_1\cn^{-1})-\frac{K_2}{2\epO^2V(1-\constt)}
\Big[1 - \frac 1{1+(1-\constt)\frac\epsilon{\epO}}\Big]
\Big[1 + \frac 1{1+(1-\constt)\frac\epsilon{\epO}}\Big]
\nonumber
\\
&\geq\epsilon(1-K_1\cn^{-1})-\frac{K_2}{\epO^2V(1-\constt)}
\Big[1 - \frac 1{1+(1-\constt)\frac\epsilon{\epO}}\Big]
\nonumber
\\
\lbeq{chi-1}
&=\epsilon\Big(1-K_1\cn^{-1}-{K_2\lambda^3}
\frac 1{1+(1-\constt)\epsilon/\epO}
\Big).
\enalign
The upper bound in \refeq{chibd-strong} is equivalent to \refeq{chi-1}.
\qed

\bigskip \noindent
{\em Proof of \refeq{omegapc-1-strong}.}
The lower bound of \refeq{omegapc-1-strong} was proved already in
\refeq{omegapc-1}.  The upper bound follows from \refeq{chibd-strong}
with $p=0$, using the lower bound of \refeq{omegapc-1} to bound
$\epsilon = \cn p_c$ in \refeq{atildef}.
\qed

\subsection{Upper bound on the largest subcritical cluster}
\label{sec-scp1}

\noindent
{\em Proof of the upper bound of
\refeq{cmaxbd1} and of \refeq{cmaxbd2}.}
We will prove that
    \eq\lbeq{cmax.1}
    \Exp\Big(|\Cmax|\Big)
    \leq
    2\chi^2(p)\log(V/\chi^3(p))
    \en
and
    \eq\lbeq{cmax.2}
    \Pro\Big(|\Cmax|\leq 2\chi^2(p)\log(V/\chi^3(p))\Big)
    \geq
    1 - \frac {\sqrt e}{[2\log(V/\chi^3(p)]^{3/2}},
    \en
if $\chi(p)\leq e^{-2}V^{1/3}$.
The desired bounds follow immediately from
\refeq{cmax.1} and \refeq{cmax.2},
provided $\lambda=V^{-1/3}\chi(p_c)\leq e^{-2}$.
However, it follows from \refeq{lambda-bd}
and our assumptions $\constt\leq b_0$ and
$V^{-1/3}\leq \lambda b_0$ that
$\lambda^3\leq b_0 +\lambda^3 b_0^3$,
which gives $\lambda^3\leq b_0(1-b_0^3)^{-1} \leq e^{-2}$.

To prove \refeq{cmax.2}, we let
$A=\log(V/\chi^3(p))$
and $k=2A\chi^2(p) $.  By assumption,
$A\geq 6\geq 1/2$, and hence $k\geq \chi^2(p)$.
We can therefore apply \refeq{cmax.3} and \refeq{C0ub} to obtain
\eq\lbeq{cmax.5}
   \Pro\Big(|\Cmax|\geq 2A\chi^2(p)\Big)
   \leq
   \frac{VP_{\geq 2A\chi^2(p)}}{2A\chi^2(p)}
   \leq \frac{V\sqrt e}{(2A)^{3/2}\chi^3(p)}
   e^{-A}=
   \frac{\sqrt e}{(2A)^{3/2}},
   \en
which is the desired bound \refeq{cmax.2}.

To prove \refeq{cmax.1}, we set
$k=2(A-1)\chi^2(p)$.  Combining
\refeq{cmax.4} and \refeq{C0ub} leads to
\eq
\begin{aligned}
   \Exp\Big(|\Cmax|\Big)
   &\leq 2(A-1)\chi^2(p)
   \Bigl[1+\frac{V\sqrt e}{(2(A-1))^{3/2} \chi^3(p)}
   e^{-A+1}\Bigr]
   \\
   &= 2Y\chi^2(p)\log(V/\chi^3(p))
\end{aligned}
\en
with
\eq Y=\Bigl(1-\frac 1A\Bigr)
\Bigl[1+\Bigl(\frac e{2(A-1)}\Bigr)^{3/2}\Bigr].
\en
To complete the proof, it suffices to show that
\eq\lbeq{cmax.6}
\bigg(\frac e{2(A-1)}\bigg)^{3/2}
\leq \frac 1A,
\en
since this implies that $Y\leq 1$.
To prove \refeq{cmax.6}, we use the
monotonicity of the function $x\mapsto (x-1)^3/x^2$
and the fact that $A\geq 6$ to conclude that
\eq
\frac {8(A-1)^3}{A^2}\geq \frac{1000}{36}
\geq e^3.
\en
\qed

\section{The magnetization}
\label{sec-M}

In this section, we prove Theorem~\ref{main-thm-mag}.  This theorem
provides upper and lower bounds on the
magnetization, which is defined by
    \eq
    \lbeq{Mdef}
    M(p,\gamma)
    =\sum_{k=1}^{V} [1-(1-\gamma)^{k}] \mathbb P_{p}(|C(0)|=k).
    \en
For fixed $p$, the function $M(p,\cdot)$ is strictly
increasing, with $M(p,0)=0$ and $M(p,1)=1$.  We denote the inverse
function by $\gamma(m)$, so that $M(p,\gamma(m))=m$ for all $m
\in [0,1]$. In addition, for $\gamma \in (0,1)$,
$M(p,\gamma)$ is strictly increasing in $p$.
Finally, recalling \refeq{chiMpg},
we note that $\partial M/\partial\gamma=(1-\gamma)^{-1}\chi$ is
monotone decreasing in $\gamma$.  Since $M(p,0)=0$
this implies that
 \eq\lbeq{M-gamma-chi}
\frac\gamma{1-\gamma}\chi(p,\gamma)
\leq M(p,\gamma)
\leq\gamma\chi(p,0).
\en

\subsection{Bounds on the magnetization}
\label{sec-mag}

We formulate our
results in the general setting of a connected transitive graph $\gr$ with
$V$ vertices and degree $\cn$,
not necessarily obeying the triangle condition \refeq{tria-con}.
Instead, we will assume that one or several of the following conditions
hold:
\eqalign
\lbeq{small-pc}
p_c &\leq \constp,
\\
\lbeq{opc-1-up}
\cn p_c & \leq 1+\constopup,
\\
\lbeq{opc-1-low}
\cn p_c & \geq 1-\constoplow,
\enalign
and last  but not least,
the triangle condition \refeq{tria-con} itself.
The constants
$\constt$, $\constp$,
$\constopup$ and $\constoplow$  in the following statements
refer to these assumptions, and when a constant is not mentioned
in a theorem, the corresponding assumption is not used.

Note that when we do assume the triangle condition, then the assumptions
\refeq{small-pc}--\refeq{opc-1-low} all follow, provided
$V$ is large enough.
To see this, we note that for any bond $\{x,y\}\in \mathbb B$,
we have $p \leq \tau_p(x,y) \leq \nabla_p(x,y)$ (just take
$u=v=x$ in \refeq{tria-def}), and hence
\eq\lbeq{pc-less-constt}
p_c\leq \constt
\en
whenever the triangle condition \refeq{tria-con} holds.
In addition, \refeq{opc-1-up}--\refeq{opc-1-low} follow from
\refeq{omegapc-1}.
Therefore,  in particular, the constants $\constp$,
$\constopup$ and $\constoplow$ can be made as small as desired
by assuming that $\constt$ and $\epO=\lambda^{-1}V^{-1/3}$ are
sufficiently small (as assumed in the theorems in
Section~\ref{sec-mainresults}).

The following propositions and corollaries immediately imply
Theorem~\ref{main-thm-mag}.
The first pair gives lower bounds on the magnetization, and the
second pair gives upper bounds.
For Corollary~\ref{Mgeq-cor}, we recall that
$N_\alpha = \epsilon^{-2}(\epsilon V^{1/3})^\alpha$ was defined in
\refeq{N_a}.

\begin{prop}
    \label{Mgeq}
    (i) Let $0<p<1$ and $0<\gamma <1$, and let
    $K=1+\frac{\cn p}{1-p}$.  Then
        \eq
        \lbeq{Mgeq0}
        M(p,\gamma)\geq
        \frac 1{2K}\left[
        \sqrt{4K\gamma+\chi^{-2}(p)}
        -\chi^{-1}(p)
        \right],
        \en
   so that in particular
    \eq
    \lbeq{Mgeq1}
    M(p,\gamma)\geq
    \frac{\sqrt{4K+1}-1}{2K}
    \min\{\sqrt\gamma,\gamma\chi(p)\}.
    \en
\\
    (ii) If $0<p_0 \leq p<1$, $0<\gamma_0  \leq \gamma <1$ and
    $0<\tilde\alpha<1$, then
         \eq\lbeq{Mgeq5}
         M(p,\gamma)\geq
         \min\bigg\{
         \Big(\frac\gamma{\gamma_0}\Big)^{\tilde\alpha}
         M(p_0,\gamma_0),
         \frac {p_0}p M(p_0,\gamma_0)+
         (1-\tilde\alpha)\frac{p-p_0}p
         \bigg\}.
         \en
\end{prop}

\begin{cor}
\label{Mgeq-cor}
Assume that $\constp$ and $\constopup$ are sufficiently small.

\noindent
i) If $0\leq \gamma\leq 1$ and $p\leq p_c$, then
        \eq
        \lbeq{Mgeq2}
        M(p,\gamma)\geq
        \frac 13
    \min\{\sqrt\gamma,\gamma\chi(p)\}.
    \en

\noindent
ii) Let  $0\leq\alpha<1$, $\tilde\alpha=(2-\alpha)^{-1}$, $\rho>0$,
and $p=p_c+\cn^{-1}\epsilon$.
Let $b_{13} = \lambda^{-1-\alpha \tilde \alpha}\rho^{-\tilde \alpha}$.
If $b_{13}V^{-1/3} \leq \epsilon \leq 1$ then
    \eq
    \lbeq{Mgeq7}
    M(p,\rho N_\alpha^{-1})
    \geq
        \frac {\epsilon}{3}
    \min\{(1-\tilde\alpha),\rho^{\tilde\alpha}
    \lambda^{\alpha \tilde \alpha}\}    .
    \en
\end{cor}

\begin{lemma}
\label{lem-Mleq}
If $\constt$ and $\constoplow$ are sufficiently
small, $p\leq p_c$ and $0\leq\gamma\leq 1$, then
    \eq
    \lbeq{Mleq2}
    M(p,\gamma)\leq
    \min\{\sqrt{12\gamma},\,\gamma\chi(p)\}.
    \en
\end{lemma}

\begin{prop}
\label{Mleq}
If $\constt$ and $\constoplow$ are sufficiently
small, $p=p_c+\cn^{-1}\epsilon\geq p_c$ and $0\leq\gamma\leq 1$, then
    \eq
    \lbeq{Mleq1}
    M(p,\gamma)\leq
    \sqrt{12\gamma}+13\epsilon.
    \en
\end{prop}

\smallskip \noindent
{\em Proof of Theorem~\ref{main-thm-mag}.}
This is an immediate consequence of Corollary~\ref{Mgeq-cor},
Lemma~\ref{lem-Mleq} and Proposition~\ref{Mleq}.
\qed

Note that for $p\leq p_c$ the lower bound \refeq{Mgeq2} and
the upper bound \refeq{Mleq2} differ only by a constant,
for all $\gamma$.
For $p\geq p_c$, our results are much weaker:
If we specialize to $\gamma$ proportional to $N_\alpha^{-1}$,
and assume that $\epsilon V^{1/3}$ is large enough (in particular, this
implies that $N_\alpha^{-1} \leq \epsilon^2$),
then our lower and upper bounds \refeq{Mgeq7} and \refeq{Mleq1}
match.

Our bounds on the magnetization are proved using the
three differential
inequalities stated in the next lemma.

    \begin{lemma}
    \label{ineq}
    If $0< p < 1$ and $0<\gamma <1$, then
        \eq
        \label{ineq1}
        (1-p)\frac{\partial M}{\partial p}
        \leq \cn (1-\gamma) M \frac{\partial M}{\partial \gamma},
        \en
        \eq
        \label{ineq2}
        M
        \leq \gamma \frac{\partial M}{\partial \gamma}
        +M^2 + pM\frac{\partial M}{\partial p},
        \en
    and
        \eq
        \lbeq{rdi}
         M \geq
        \left[{\cn \choose 2}p^2(1-p)^{\cn -2}(1-\nabla_p^{\rm max})^3
        -p-\nabla_p^{\rm max} \right]p\cn(1-\gamma)
    M^2 \frac{\partial M}{\partial\gamma},
        \en
    where $\nabla_p^{\rm max} = \max_{x,y \in \mathbb V}\nabla_p (x,y)$.
    \end{lemma}

The differential inequalities \eqref{ineq1}--\eqref{ineq2}
were derived and used by Aizenman and Barsky
\cite{AB87} to prove sharpness of the percolation phase transition
on $\Z^n$, and will be used to prove our lower bounds
on $M(p,\gamma)$.  The derivations in \cite{AB87} extend without difficulty
to an arbitrary transitive graph.
The differential inequality \refeq{rdi}, which is
a
variant
of an inequality derived by Barsky and Aizenman \cite{BA91}, will
be used to prove our upper bounds on $M(p,\gamma)$.
We give a proof of \refeq{rdi} in Appendix~\ref{sec-M-diff-ine}.

\subsection{Lower bounds on the magnetization}
\label{sec-Mlb}

In this section, we prove Proposition \ref{Mgeq}
and Corollary~\ref{Mgeq-cor}, using the first two differential
inequalities of Lemma~\ref{ineq}.

{\em\noindent Proof of Proposition \ref{Mgeq}.}
 (i)
We fix
$p\in(0,1)$, and drop the $p$ dependence from the notation.
Inserting \eqref{ineq1}
into \eqref{ineq2}, and using $\tilde K=\frac{\cn p}{1-p}$
and $1-\gamma \leq 1$, we get
    \eq
    \begin{aligned}
    \label{ineq3}
    M
    &\leq
    \gamma\frac {dM}{d\gamma} +M^2+
    \tilde KM^2\frac{dM}{d\gamma}.
    \end{aligned}
    \en
Since $M>0$ as long as $\gamma>0$, we get
    \eq
    \frac 1M \frac{d\gamma}{dM}
    -\frac 1{M^2} \gamma
    \leq
    \tilde K
    +
    \frac{d\gamma}{dM},
    \en
where we are using the fact that $M$ has a well-defined inverse function.
Therefore,
    \eq
    \label{inv-ineq1}
    \frac d{dM}
    \left(\frac \gamma M \right)
    \leq
    \tilde K
    +
    \frac{d\gamma}{dM}.
    \en

Next we integrate \eqref{inv-ineq1} and use
that $\gamma(0)=0$ and $\lim_{M\to 0}\frac {\gamma(M)}{M}
=\gamma^\prime(0)= 1/M'(0)=\chi^{-1}(p)$ to get
    \eq
    \label{int-ineq1a}
    \frac \gamma M \leq \chi^{-1}+\tilde KM + \gamma
    \en
where we used the shorthand $\chi^{-1}$ for
$\chi^{-1}(p)$.
Observing that $1-(1-\gamma)^k\geq 1-(1-\gamma)=\gamma$,
we see from \refeq{Mdef} that $\gamma\leq M$,
which simplifies \eqref{int-ineq1a}
to
\eq
    \label{int-ineq1}
    \frac \gamma M \leq \chi^{-1}+KM,
    \en
where $K=\tilde K+1$.
Multiplying by $M/K$ and completing
the square on the right side, we thus obtain
    \eq
    \label{int-ineq2}
    \frac \gamma {K}
    +
    \left[\frac{\chi^{-1}}{2K}\right]^2
    \leq
    \left[
    M+\frac{\chi^{-1}}{2K}
    \right]^2 .
    \en
Since $M \geq 0$, this implies that
    \eq
    \label{alt2}
    M
    \geq
    \sqrt{\frac \gamma{K}
    +
    \left[\frac{\chi^{-1}}{2K}\right]^2}
    -\frac{\chi^{-1}}{2K}.
    \en
This completes the proof of \refeq{Mgeq0}.

To prove \refeq{Mgeq1}, let us first assume that
$\gamma\geq \chi^{-2}(p)$.  By \refeq{Mgeq0} and the fact that the function
$f(x)=\frac 1{\sqrt x}(\sqrt {x+\chi^{-2}}-\chi^{-1})$ is increasing,
we conclude that
    \eqalign
    M(p,\gamma)
    & \geq
    \sqrt{\frac \gamma K} f(4K\gamma)
    \nnb
    & \geq \sqrt{\frac \gamma K} f(4K\chi^{-2}(p))
    =
    \frac{\sqrt{4K+1}-1}{2K}
    \sqrt\gamma
    \nnb
    & =
    \frac{\sqrt{4K+1}-1}{2K}
    \min\{\sqrt\gamma,\gamma\chi(p)\}.
    \enalign
On the other hand, if $\gamma\leq \chi^{-2}(p)$, we use the fact that
the function
$g(x)= \frac 1{ x}(\sqrt {x+\chi^{-2}}-\chi^{-1})$ is decreasing,
together with the bound $M(p,\gamma) \geq 2\gamma g(4K\gamma)$ of
\refeq{Mgeq0},
to arrive at the same conclusion.
This completes the proof of (i).

\medskip \noindent
(ii)
The result is immediate if $\gamma_0=\gamma$ or $p_0=p$ so we
assume that $\gamma_0 < \gamma$ and $p_0<p$.
We rewrite \eqref{ineq2} as
    \eq
    \label{ineq2a}
    0
    \leq
    \frac 1M\frac{\partial M}{\partial\gamma}
    +\frac 1\gamma\frac{\partial }{\partial p}(pM-p),
    \en
and then integrate \eqref{ineq2a} over the rectangle
$[\gamma_0,\gamma]\times [p_0,p]$.  This yields
    \eq
    \label{partially-integrated}
    0\leq
    \int^p_{p_0}d\tilde p\log
    \left(
    \frac{M(\tilde p,\gamma)}{M(\tilde p,\gamma_0)}
    \right)+
    \int_{\gamma_0}^\gamma d\tilde\gamma\frac 1{\tilde\gamma}
    \left(pM(p,\tilde\gamma)-p_0M(p_0,\tilde\gamma)
    -(p-p_0)\right).
    \en
Since
    \eq
    \label{monoton}
    0\leq M(p_0,\gamma_0)
    \leq M(\tilde p,\tilde\gamma)
    \leq M(p,\gamma)
    \en
whenever
$(\tilde p,\tilde\gamma)\in[\gamma_0,\gamma]\times [p_0,p]$,
it follows that
    \eq
    \label{integrated}
    0\leq
    (p-{p_0})\log
    \left(
    \frac{M(p,\gamma)}{M(p_0,\gamma_0)}
    \right)+
    \log\left(\frac\gamma{\gamma_0}\right)
    \left(pM(p,\gamma)-p_0M(p_0,\gamma_0)
    -(p-p_0)\right).
    \en
Dividing by $\log(\gamma/{\gamma_0})$,
we conclude that
    \eq
    \lbeq{Mgeq3}
    M(p,\gamma)
    \geq
\frac {p_0}pM(p_0,\gamma_0)+
\frac{p-p_0}{p}
    \left[1-
    \frac{\log\{M(p,\gamma)/M(p_0,\gamma_0)\}}%
    {\log\{\gamma/\gamma_0\}}
    \right].
    \en

If $M(p,\gamma)/M(p_0,\gamma_0)
\leq (\gamma/\gamma_0)^{\tilde\alpha}$, then \refeq{Mgeq3} gives
\eq
M(p,\gamma)\geq
\frac {p_0}pM(p_0,\gamma_0)+
\frac{p-p_0}{p}\big[1-\tilde\alpha\big].
\en
If, on the other hand, $M(p,\gamma)/M(p_0,\gamma_0)
\geq (\gamma/\gamma_0)^{\tilde\alpha}$, then it is trivially the
case that
\eq
M(p,\gamma)\geq M(p_0,\gamma_0)(\gamma/\gamma_0)^{\tilde\alpha}.
\en
Therefore, as desired,
         \eq
     \lbeq{Mgeq5-rep}
         M(p,\gamma)\geq
         \min\bigg\{
         \Big(\frac\gamma{\gamma_0}\Big)^{\tilde\alpha}
         M(p_0,\gamma_0),
         \frac {p_0}p M(p_0,\gamma_0)+
         (1-\tilde\alpha)\frac{p-p_0}p
         \bigg\}.
         \en
\qed

\bigskip
{\em\noindent Proof of Corollary~\ref{Mgeq-cor}.}
(i) The function $K=K(p)=1+\cn p/(1-p)$
is increasing in $p$, so $K(p) \leq K(p_c)$
for $p \leq p_c$.  Since the function $(\sqrt{4K+1}-1)/2K$ is
decreasing in $K$, for a lower bound we can replace $K$ by $K(p_c)$
in \refeq{Mgeq1}.  Since $K(p_c) \to 2$ as $\constp$ and $\constopup$
go to zero, \refeq{Mgeq2} then follows.

\medskip \noindent
(ii) We apply Proposition~\ref{Mgeq}(ii), whose conclusion is
repeated above in \refeq{Mgeq5-rep}, with $p_0=p_c$,
$\gamma = \rho N_\alpha^{-1}$ and
$\gamma_0=\epO^2=\chi^{-2}(p_c)$.  The requirement $\gamma \geq \gamma_0$
for \refeq{Mgeq5-rep} is equivalent to our hypothesis that
$\epsilon \geq b_{13}V^{-1/3}$.  It suffices to show that
\eq
\lbeq{suff1}
\Big(\frac\gamma{\gamma_0}\Big)^{\tilde\alpha}
         M(p_c,\gamma_0)
     \geq
     \frac{\epsilon}{3}\rho^{\tilde \alpha}\lambda^{\alpha \tilde\alpha}
\en
and
\eq
\lbeq{suff2}
    \frac {p_c}p M(p_c,\gamma_0)+
         (1-\tilde\alpha)\frac{p-p_c}p
     \geq
     \frac{\epsilon}{3}(1-\tilde\alpha).
\en
For \refeq{suff1}, we use \refeq{Mgeq0} and the observation in
the proof of part (i) to see that
\eqalign
M(p_c,\gamma_0)
&\geq
\frac {\chi^{-1}(p_c)}{2K}
\left[\sqrt{4K+1}-1 \right] \geq \frac 13 \epsilon_0
\enalign
if $a_1$ and $a_2$ are sufficiently small.
Since $(\gamma/\gamma_0 )^{\tilde\alpha}=
\rho^{\tilde\alpha}\lambda^{\alpha \tilde \alpha}({\epsilon}/{\epsilon_0})$,
\refeq{suff1} follows.
For \refeq{suff2}, we bound the first term on the left side below by zero,
and note that
\eq
\lbeq{Mgeq6}
\frac{p-p_c}p
=\frac{\epsilon}{\epsilon+\cn p_c}
\geq \frac{\epsilon}3,
\en
since $\epsilon\leq 1$ by assumption and
$\cn p_c\leq 2$ if $\constopup$ is small enough.
\qed

\subsection{Upper bounds on the magnetization}
\label{sec-Mub}

We now prove Lemma~\ref{lem-Mleq} and Proposition~\ref{Mleq}.
Lemma~\ref{lem-Mleq} is proved by integration of the differential
inequality \eqref{ineq2}, assuming the triangle condition.
We then use the extrapolation principle
of \cite{AB87,AF86,BA91} to convert the upper bound on $M(p_c,\gamma)$
to an upper bound valid for $p>p_c$.  This is perhaps surprising,
since $M$ is an increasing function of $p$.  However, it is
also increasing
in $\gamma$, and we will see that it is possible to use the differential
inequality \eqref{ineq1} to compensate for an increase in $p$ with
a decrease in $\gamma$.

\medskip\noindent{\it Proof of Lemma~\ref{lem-Mleq}.}
We first note that $M(p,\gamma) \leq \gamma \chi(p)$ for all $p$ and
$\gamma$, by \refeq{M-gamma-chi}.  Since $M(\cdot, \gamma)$ is
increasing, it suffices to prove that
\eq
\lbeq{Mdeltaup}
    M(p_c, \gamma) \leq \sqrt{12\gamma}.
\en
We assume that the triangle condition is
satisfied and that \refeq{opc-1-low} holds for a sufficiently small
constant $\constoplow$. Under these conditions,
$1-\constoplow\leq\cn p_c\leq (1-\constt)^{-1}$,
$p_c\leq\constt$ and $\cn=\cn p_c/p_c\geq (1-\constoplow)/\constt$,
so \refeq{rdi}
implies that
    \eq\lbeq{Mrev}
    M(p_c,\gamma)
    \geq
    \frac 1{2e}\big[1-\bigo(\constt \vee\constoplow)\big]
    (1-\gamma)
    \frac{d M(p_c,\gamma)}{d\gamma}
    M(p_c,\gamma)^2.
\en
For $p=p_c$, this gives
\eq
\lbeq{dM2}
    \frac{1}{2} \frac{dM^2}{d\gamma} \leq \frac{2e}{1-\gamma}
    \big[1-\bigo(\constt \vee \constoplow)\big].
\en
We integrate \refeq{dM2} over the interval $[0,\gamma]$, using
$M(p_c,0)=0$, to see that
\eq
    M^2(p_c,\gamma) \leq \frac{4e \gamma}{1-\gamma}
    \big[1-\bigo(\constt \vee \constoplow)\big].
\en
For $\gamma \in [0,\frac 1{12}]$, this implies \refeq{Mdeltaup},  provided
$\constt$ and $\constoplow$ are sufficiently small.
Finally, we note that we
can remove the restriction $\gamma \in [0,\frac 1{12}]$,
since trivially, $M(p,\gamma)\leq 1\leq \sqrt{12\gamma}$
if $\gamma\geq \frac 1{12}$.
\qed

\medskip\noindent{\it Proof of Proposition~\ref{Mleq}.}
Following \cite{BA91},
we apply the extrapolation principle used in \cite{AB87},
to extend \refeq{Mdeltaup} to \refeq{Mleq1}.
The extrapolation principle is explained in \cite{AF86}.
In our setting,
the finite size effect will need to be taken into account.
We find it most convenient to use the variable $h$ rather than $\gamma$,
and define $\tilde{M}(p,h) = M(p,1-e^{-h})$, for $h \geq 0$.

Assuming that $\epsilon\leq 1$,
the differential inequality \eqref{ineq1} implies that
\eq
    \frac{\partial \tilde{M}} {\partial p}
    \leq A\cn\tilde{M} \frac{\partial \tilde{M}} {\partial h}
\en
where $A=(1-p_c-\cn^{-1})^{-1}=1+O(\constt)+O(\constoplow)$.
For fixed $m \in [0,1]$ and fixed $p\in (0,1)$, we can solve the equation
$\tilde{M}(p,h)=m$ for $h = h(p)$, so that $\tilde{M}(p,h(p))=m$.
Differentiation of this identity with respect to $p$ gives
\eq
    \frac{\partial \tilde{M}} {\partial p} +
    \frac{\partial \tilde{M}} {\partial h}
    \left. \frac{\partial h} {\partial p} \right|_{\tilde{M}=m}
     = 0.
\en
Therefore,
\eq
\lbeq{dhdp}
    0 \leq - \left. \frac{\partial h} {\partial p} \right|_{\tilde{M}=m}
    =
    \frac{\frac{\partial \tilde{M}} {\partial p}}
    {\frac{\partial \tilde{M}} {\partial h}}
    \leq A\cn m.
\en

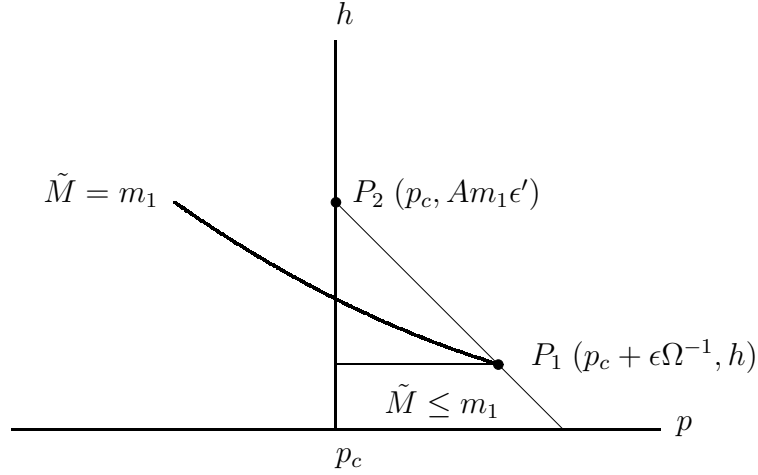
\begin{figure}
\begin{center}
\setlength{\unitlength}{0.0170in}%
\begin{picture}(200,160)(50,750)

\thicklines
\put(50,760){\line( 1, 0){ 200}}
\put(150,760){\line(0,1){120}}
\qbezier(100,830)(150,795)(200,780)
\thinlines
\put(150,830){\line(1,-1){70}}
\put(150,780){\line(1,0){50}}
\put(150,830)
{\makebox(0,0)[lb]{\raisebox{0pt}[0pt][0pt]{\circle*{3}}}}
\put(200,780)
{\makebox(0,0)[lb]{\raisebox{0pt}[0pt][0pt]{\circle*{3}}}}
\put(255,760)
{\makebox(0,0)[lb]{\raisebox{0pt}[0pt][0pt]{$p$}}}
\put(150,750)
{\makebox(0,0)[lb]{\raisebox{0pt}[0pt][0pt]{$p_c$}}}
\put(210,780)
{\makebox(0,0)[lb]{\raisebox{0pt}[0pt][0pt]{$P_1\;(p_c+\epsilon
\cn^{-1},h)$}}}
\put(155,830)
{\makebox(0,0)[lb]{\raisebox{0pt}[0pt][0pt]{$P_2\;(p_c,Am_1\epsilon')$}}}
\put(150,885)
{\makebox(0,0)[lb]{\raisebox{0pt}[0pt][0pt]{$h$}}}
\put(165,765)
{\makebox(0,0)[lb]{\raisebox{0pt}[0pt][0pt]{$\tilde{M}\leq m_1$}}}
\put(60,830)
{\makebox(0,0)[lb]{\raisebox{0pt}[0pt][0pt]{$\tilde{M}= m_1$}}}

\end{picture}
\end{center}
\caption{The extrapolation geometry.}
\label{fig-extrapolation}
\end{figure}

The upper and lower bounds of \refeq{dhdp} imply that a contour
line $\tilde{M}=m_1$ in the $(p,h)$-plane (with $p$-axis
horizontal and $h$-axis vertical) passing through a point
$P_1=(p_1,h_1)$ is such that $\tilde{M}(P) \leq m_1$ for all
points $P$ in the first quadrant that lie on or below the lines of
slope 0 and $-A\cn m_1$ through $P_1$; see
Figure~\ref{fig-extrapolation}. In addition, if $P_2=(p_2,h_2)$ is
on the line through $P_1$ with slope $-A\cn m_1$, with $p_2<p_1$,
and if we set $m_2 = \tilde{M}(P_2)$, then $m_2 \geq m_1$. For if
to the contrary $m_2<m_1$, then $P_1$ would lie below the line
through $P_2$ of slope $-A\cn m_2$, which would imply that $m_1
\leq m_2$, a contradiction. We will use the fact that $m_2 \geq
m_1$ below.

Fix $h$, and fix $\epsilon >0$.
Let $P_1=(p_c+\epsilon \cn^{-1},h)$, and define
$m_1 = m_1(\epsilon)= \tilde{M}(P_1)$.
Let
\eq
    \epsilon' = \epsilon + \frac{h}{Am_1},
\en
and define $P_2 = (p_c,A m_1 \epsilon')$ and $m_2=\tilde{M}(P_2)$.
The points $P_1$ and
$P_2$ are collinear on the line through $P_1$ with slope
$-A\cn m_1$.  Therefore, as observed above, $m_1 \leq m_2$.
Applying \refeq{Mdeltaup} gives
\eqalign
    \tilde{M}(p_c+\epsilon \cn^{-1}, h)
    &=m_1
    \leq m_2
    =  \tilde{M}(p_c, A m_1 \epsilon')
    \leq \sqrt{12} (1-e^{-A m_1 \epsilon'})^{1/2}
    \nonumber \\
    &=
    \sqrt{12}
    \big(1-e^{-Am_1 \epsilon}+e^{-Am_1 \epsilon}[1-e^{-h}]\big)^{1/2}
    \nonumber \\
    &\leq \sqrt{12} (Am_1 \epsilon+\gamma)^{1/2}
    ,
\enalign
with $\gamma=1-e^{-h}$.  The inequality
\eq
    m_1^2 \leq 12(Am_1\epsilon + \gamma)
\en
has roots
\eq
    m^{\pm} =  6A\epsilon \pm \sqrt{12\gamma+  (6A\epsilon)^2}.
\en
The root $m^+$ is positive and $m^-$ is negative.
Thus we have
\eq
\lbeq{Mlower1}
    M(p_c+\epsilon \cn^{-1}, \gamma)
    = m_1
    \leq m^
    + \leq 6A\epsilon + \sqrt{12\gamma+ (6A\epsilon)^2}
    \leq 12A\epsilon + \sqrt{12\gamma}.
\en
This completes the proof of \refeq{Mleq1}, since we can
choose $A$ arbitrarily close to $1$ by choosing $\constt$ and
$\constoplow$ sufficiently small.
\qed

\section{The cluster size distribution}
\label{sec-proof}

In this section, we prove
Theorems~\ref{main-thm-critical}~i), \ref{main-thm-sup}~i)
and \ref{main-thm-pp}~i).
The magnetization
    \eq
    M(p,\gamma)=
    \sum_{k=1}^{V} [1-(1-\gamma)^{k}] \mathbb P_{p}(|C(0)|=k)
    \en
is a generating function for $\Pbold_p(|C(0)|=k)$.  In the spirit
of a Tauberian theorem, we will
use the bounds on $M(p,\gamma)$
established in Section~\ref{sec-M} to obtain bounds on
$P_{\geq s}(p)$.  We recall the upper bound
 \eq \lbeq{Mub}
    M(p,\gamma) \leq
    \sqrt{12\gamma} + 13\epsilon ,
    \en
proved in \refeq{Mleq1} for all $p\geq p_c$ provided
$\constt$ and $\constoplow$ are sufficiently small,
and the lower bound
     \eq
     \lbeq{Mgeq2-rep}
     M(p,\gamma)\geq
     \frac 1{3}\min\{\sqrt\gamma,\gamma\chi(p)\},
     \en
proved in \refeq{Mgeq2} for all $p\leq p_c$ provided $\constp$ and
$\constopup$ are sufficiently small.
The discussion of Section~\ref{sec-mag} shows that the constants
$\constp$, $\constopup$ and $\constoplow$ can be made
arbitrarily small if $\constt\leq b_0$, $\epO\leq b_0$
and $b_0$ is chosen small enough.

The cluster size distribution and magnetization are related by the
following lemma.

\begin{lemma}
\label{lem-csd-M}
Let $p \in [0,1]$, $k >0$ and $0 \leq \gamma, \tilde\gamma \leq 1$.
Then
    \eqalign \lbeq{PleqM}
    P_{\geq k}(p)
    &\leq \frac e{e-1}M(p,k^{-1}),
    \\
    \lbeq{PgeqM}
    P_{\geq k }(p)
    & \geq
    M(p,\gamma)
    -  \frac \gamma{\tilde \gamma}\,e^{\tilde \gamma k}
    M(p,\tilde \gamma).
    \enalign
\end{lemma}

\proof
The bound \refeq{PleqM}
follows immediately from the definition of $M$ and the fact that
$1-e^{-1}\leq 1-(1-k^{-1})^\ell$ whenever $\ell\geq k$.

To prove \refeq{PgeqM}, we note that
$[1-(1-\gamma)^{\ell }] \leq \ell \gamma$.  Also,
$\ell \tilde \gamma \leq e^{\ell
\tilde \gamma}-1 =e^{\ell \tilde \gamma}(1-e^{-\ell \tilde
\gamma})$, which combined with $e^{-\tilde\gamma}\geq
1-\tilde\gamma$ gives
$\ell\tilde\gamma\leq{e^{\ell\tilde\gamma}}(1-(1-\tilde
\gamma)^\ell)$. Therefore,
    \eqalign
    M(p,\gamma)
    &=\sum_{\ell=1}^{V} (1-(1-\gamma)^\ell )
    \,{\mathbb P}_{p}(|C(0)|=\ell )
    \nonumber \\
    &\leq \gamma \sum_{\ell \leq k} \ell
    \,{\mathbb P}_{p}(|C(0)|=\ell )
    + \sum_{\ell \geq k}  \,{\mathbb P}_{p}(|C(0)|=\ell )
    \nonumber \\
    &\leq \frac \gamma{\tilde \gamma}\,e^{\tilde \gamma k}
    M(p,\tilde \gamma)
        + P_{\geq k }(p).
    \enalign
\qed

We will use
\refeq{Mub}--\refeq{Mgeq2-rep}
and Lemma~\ref{lem-csd-M}
to prove the bounds in the following lemma.

\begin{lemma}\label{lem-Pgeqk}
There is a constant $b_0$ such that
the following statements hold
provided $\epO\leq b_0$ and the
triangle condition \refeq{tria-con}
is valid with $\constt\leq b_0$.

\noindent
i) If $p=p_c+\cn^{-1}\epsilon\geq p_c$ then
\eq
\lbeq{Pgekbd-up1}
P_{\geq k}(p)\leq
21 \epsilon + 6k^{-1/2}.
\en

\noindent
ii) If $p\leq p_c$ then
    \eq
    \lbeq{Pgekbd}
    \frac 1{360}k^{-1/2}
    \leq P_{\geq k}(p)\leq
    6k^{-1/2}.
    \en
provided $k\leq \frac 1{3600}\chi(p)^2$ for the lower bound (this assumption
is not needed for the upper bound).

\noindent
iii) If $p=p_c+\epsilon\cn^{-1}$ and
$k\leq [100(|\epsilon|+\epsilon_0)]^{-2}$
then
\eq
    \lbeq{Pgekbd2}
    \frac 1{360}k^{-1/2}
    \leq P_{\geq k}(p)\leq
    6k^{-1/2}.
    \en

\end{lemma}

\proof
(i)
Inserting \refeq{Mub} into \refeq{PleqM} gives \refeq{Pgekbd-up1}.

\noindent (ii)
For the upper bound in \refeq{Pgekbd}, we take $p \leq p_c$
and note that $P_k(p) \leq P_k(p_c) \leq 6k^{-1/2}$, using monotonicity
in the first step and \refeq{Pgekbd-up1} in the second.
For the lower bound, we
apply \refeq{PgeqM} with $\tilde\gamma=1/k$.  Since
\eq
M(p,k^{-1})\leq M(p_c,k^{-1})\leq\sqrt{12/ k}
\en
by \refeq{Mub}, \refeq{PgeqM} implies that
\eq
    P_{\geq k}(p) \geq M(p,\gamma) -
    \sqrt{12 k} \gamma e  .
\en
If $\gamma\geq\chi^{-2}(p)$, then \refeq{Mgeq2-rep} implies
$M(p,\gamma) \geq \frac 13 \sqrt{\gamma}$, and hence
\eq
P_{\geq k}(p)\geq
 \frac13{\sqrt \gamma} - \sqrt{12 k}\gamma e
\geq
 \frac13{\sqrt \gamma} \big(1-30 \sqrt{\gamma k}\big).
\en
The choice $\gamma=\frac 1{60^2 k}$ gives the lower bound
of \refeq{Pgekbd}.

\noindent (iii)
To prove the upper bound of \refeq{Pgekbd2}, we note that
$k\leq [100(|\epsilon|+\epsilon_0)]^{-2}$ implies
$|\epsilon| \leq \frac 1{100}k^{-1/2}$.
It then follows from \refeq{PleqM} and \refeq{Mub} that
\eq
    P_{\geq k}(p) \leq P_{\geq k}(p_c+ 0.01 k^{-1/2}\cn^{-1})
    \leq \frac{e}{e-1}
    \left[
    \sqrt{12} + 0.13
    \right] k^{-1/2},
\en
which gives the desired bound.
For the lower bound, we note that
$P_{\geq k}(p) \geq P_{\geq k}(p_c-|\epsilon|\cn^{-1})$
and that the condition on $k$ in \refeq{Pgekbd2} implies the
condition on $k$ in \refeq{Pgekbd} for $p_c-|\epsilon|\cn^{-1}$,
by the lower bound on the susceptibility in \refeq{chibd}.
Therefore \refeq{Pgekbd2} follows from \refeq{Pgekbd}.
\qed

\bigskip \noindent
{\em Proof of Theorem~\ref{main-thm-critical} i).}
Lemma~\ref{lem-Pgeqk}~iii) immediately implies the statement of
Theorem~\ref{main-thm-critical}~i) with
$b_1=[100(\Lambda+\lambda^{-1})]^{-2}$, $b_2=1/360$ and $b_3=6$.
\qed

\smallskip \noindent
{\em Proof of Theorem~\ref{main-thm-sup} i).}
We set $k=V^{2/3}$ in
\refeq{cmax.4} and apply \refeq{Pgekbd-up1} to obtain, as required,
    \eq\lbeq{cmax.1A}
    \Exp\Big(|\Cmax|\Big)
    \leq
    21\epsilon V+7V^{2/3}.
    \en
The bound \refeq{cmax.2A} then follows from Markov's inequality.
\qed


Recall that
the percolation probability $\theta_\alpha$
is defined, for $p> p_c$ and $0 < \alpha < 1$, by
\eq\lbeq{theta_a-gen}
\theta_\alpha(p)=P_{\geq N_\alpha}(p),
\en
with $N_\alpha= \epsilon^{-2}(\epsilon V^{1/3} )^\alpha$.
We now prove
the bound \refeq{bound1} on $\theta_\alpha(p)$ of
Theorem~\ref{main-thm-pp}~i).

%
%
%

\medskip \noindent
{\em Proof of Theorem~\ref{main-thm-pp}~i).}
({\em Upper bound on $\theta_\alpha(p)$.})
If $\alpha > 0$ and $\epsilon V^{1/3}\geq 1$, then
$N_\alpha\geq\epsilon^{-2}$ and the upper bound of \refeq{bound1}
follows from \refeq{Pgekbd-up1}.

\smallskip \noindent
({\em Lower bound on $\theta_\alpha(p)$.})
 We use \refeq{PgeqM} with $k = N_\alpha$,
$\tilde\gamma = N_\alpha^{-1}$, and  $\gamma = \rho N_\alpha^{-1}$
(with $\rho >0$ to be chosen below) to obtain
\eq
    \theta_\alpha(p )
    \geq
    M(p,\rho N_\alpha^{-1})
    - \rho  e
    M(p,N_\alpha^{-1}).
\en
Let $\tilde\alpha=(2-\alpha)^{-1}$.
Let $b_{9} = \lambda^{-1-\alpha \tilde \alpha}\rho^{-\tilde \alpha}$,
and assume that $\epsilon \geq b_{9} V^{-1/3}$.
By Corollary~\ref{Mgeq-cor}~ii),
    \eq \lbeq{Mgeq7-rep}
    M(p,\rho N_\alpha^{-1})
    \geq\frac {\epsilon}{3}
    \min\{(1-\tilde\alpha),\rho^{\tilde\alpha}
    \lambda^{\alpha \tilde \alpha}\}.
    \en
Assuming that $N_\alpha \geq \epsilon^{-2}$, which follows if we
also assume $b_{9} \geq 1$, it follows from
Proposition~\ref{Mleq} that
\eq
    \rho  e M(p,N_\alpha^{-1})
    \leq \rho e (\sqrt{12}+13)\epsilon.
\en
Therefore,
\eq
    \theta_\alpha(p)\geq
    \epsilon
    \Big( \frac 13\min\{(1-\tilde\alpha),\rho^{\tilde\alpha}
    \lambda^{\alpha \tilde \alpha} \}
    -
    \rho e (\sqrt{12} + 13)\Big).
\en
Since $\alpha<1$ implies $\tilde\alpha<1$, we can make the ratio of
the first to the second term as large as we want
by choosing $\rho$ sufficiently small depending on $\alpha$ and $\lambda$.
This gives the lower bound of \refeq{bound1}, with $b_{9}$
and $b_{10}$ depending on $\alpha$ and $\lambda$.
\qed

\section{Lower bound on the largest subcritical cluster}
\label{sec-lbsc}

In this section, we complete the proof of Theorem~\ref{main-thm-sub}~ii)
by proving the lower bound of \refeq{cmaxbd1}, and \refeq{cmaxbdom}.
Given $s>0$, let $Z_{\geq s}$ be the number
of vertices that lie in clusters of size $s$ or larger,
as defined in \refeq{Zdef}.
We will use the bound on the variance of
$Z_{\geq s}$ given in the following lemma.

\begin{lemma}
\label{var-lemma}
Let $s>0$ and $p\in (0,1)$.
Then
\eqalign
\Var\left[Z_{\geq s}\right]
\leq
V\chi(p)
\,.
\lbeq{Z-var-bound}
\enalign
\end{lemma}

\begin{proof}
Let
\eq\lbeq{chigeqs}
\chi_{\geq s}(p)=\Exp\left[|C(0)|I[|C(0)|\geq s\right].
\en
We will prove that
\eq
\lbeq{Z-var-bound0}
\Exp\big(Z^2_{\geq s}\big)
\leq
\left(\Exp\left[Z_{\geq s}\right]\right)^2
+
V\chi_{\geq s}(p)\big(1-P_{\geq s}(p)\big),
\en
which implies \refeq{Z-var-bound}.

We start by rewriting the expectation of $Z^2_{\geq s}$ as
\eq
\Exp\left[ Z^2_{\geq s}\right] =
\sum_{x,y\in\mathbb V}
\sum_{S:x\in S,\atop |S|\geq s}
\sum_{T:y\in T,\atop |T|\geq s}
\Pro\big(C(x)=S,C(y)=T\big).
\en
Next, we observe that $C(x)$ and $C(y)$ must be identical if they
are not disjoint.  As a consequence, the sum decomposes into
two terms: the term
\eqalign
\lbeq{firstterm}
\sum_{x,y\in\mathbb V}
\sum_{S:x,y\in S,\atop |S|\geq s}
\Pro\big(C(x)=S\big)
=
\sum_{x\in\mathbb V}
\sum_{S:x\in S,\atop |S|\geq s}
|S|\Pro\big(C(x)=S\big)
=
V\chi_{\geq s}(p)
\enalign
and the term
\eqalign
\lbeq{secondterm}
\sum_{x\in\mathbb V}
&\sum_{S:x\in S,\atop |S|\geq s}
\sum_{y\in\mathbb V\setminus S}
\sum_{T:y\in T,\atop |T|\geq s}
\Pro\big(C(x)=S,C(y)=T\big)
\nonumber
\\
&=\sum_{x\in\mathbb V}
\sum_{S:x\in S,\atop |S|\geq s}
\Pro\big(C(x)=S\big)
\sum_{y\in\mathbb V\setminus S}
\Pro\big(|C(y)|\geq s\mid C(x)=S\big).
\enalign
Denoting the set of all edges which either join two points
in $S$ or join a point in $S$ to a point in $\mathbb V\setminus S$
by $B_+(S)$, we now rewrite the conditional probability as
\eq
\lbeq{secondcond}
\Pro\big(|C(y)|\geq s\mid C(x)=S\big)
=\Pro\big(|C(y)|\geq s\mid
\text{all edges in $B_+(S)$ are vacant}\big).
\en
By the FKG inequality, \refeq{secondcond} is bounded above by
the unconditioned probability $\Pro\big(|C(y)|\geq s\big)$.
Therefore, \refeq{secondterm} is bounded by
\eqalign
\sum_{x\in\mathbb V}
&\sum_{S:x\in S,\atop |S|\geq s}
\Pro\big(C(x)=S\big)
\sum_{y\in\mathbb V\setminus S}
\Pro\big(|C(y)|\geq s\big)
\nonumber
\\
&=V
\sum_{S:0\in S,\atop |S|\geq s}
\big(V-|S|\big)\Pro\big(C(0)=S\big)
\Pro\big(|C(0)|\geq s\big)
\nonumber
\\
&=
\left(\Exp\left[Z_{\geq s}\right]\right)^2
-V\chi_{\geq s}(p)P_{\geq s}(p)
.
\lbeq{secondtermbound}
\enalign
The combination of \refeq{firstterm} and
\refeq{secondtermbound} proves
\refeq{Z-var-bound0} and hence \refeq{Z-var-bound}.
%
\end{proof}

\smallskip \noindent
{\em Proof of \refeq{cmaxbdom}.}
Let $p\leq p_c$ and $\omega \geq 1$.
Assume that $\epO\leq b_0$ and that \refeq{tria-con}
holds for some $\constt\leq b_0$ with $b_0$ as in Lemma~\ref{lem-Pgeqk}.
We must show that
\eq
\lbeq{Cmaxlbd1}
    \Pro\Big(|\Cmax|\geq \frac{\chi^2(p)}{3600\omega}\Big)
    \geq\big(1+\frac {36\chi^3(p)}{\omega V}\Big)^{-1}.
\en

By definition,  $|\Cmax|\geq s$
if and only if $Z_{\geq s}>0$.  By the Cauchy--Schwarz inequality,
\eq
\Exp\big[Z_{\geq s}\big]
=\Exp\big[Z_{\geq s}I[Z_{\geq s}>0]\big]
\leq \sqrt{\Exp\big[Z^2_{\geq s}\big]\Pro\big(Z_{\geq s}>0\big)}
\en
and thus
\eq
\Pro\big(|\Cmax|\geq s\big)
=
\Pro\big(Z_{\geq s}>0\big)
\geq
\frac{\big(\Exp\big[Z_{\geq s}\big]\big)^2}{\Exp\big[Z^2_{\geq s}\big]}
=(1+x)^{-1}
\en
where
    \eq
    x=
    \frac{\Var\big[Z_{\geq s}\big]}{\big(\Exp\big[Z_{\geq s}\big]\big)^2}
    .
    \en
By Lemma~\ref{var-lemma}, the variance of $Z_{\geq s}$
is bounded by $V\chi(p)$.  Combined with \refeq{ZExp}, this gives
$x\leq {\chi(p)}V^{-1}[P_{\geq s}(p)]^{-2}$ and thus
    \eq
    \lbeq{Cmaxlbd2}
    \Pro\big(|\Cmax|\geq s\big)
    \geq
    \Bigl(1+\frac{\chi(p)}{V[P_{\geq s}(p)]^2 }\Bigr)^{-1}.
    \en
To complete the proof, we note that \refeq{Cmaxlbd1} is trivial if
$\omega \geq \chi^2(p)/3600$.  For $\omega \leq \chi^2(p)/3600$,
we chose $s=\chi^2(p)/3600\omega$ and use \refeq{Pgekbd} to bound
$P_{\geq s}(p)$ from below by $\frac {\sqrt\omega}{6\chi(p)}$.
This gives \refeq{Cmaxlbd1}. \qed

\smallskip \noindent
{\em Proof of the lower bound of \refeq{cmaxbd1}.}
We recall
from \refeq{lambda-bd}
that $\chi^3(p)/V\leq \chi^3(p_c)/V=\lambda^3\leq \constt+V^{-1}$.
Therefore, by \refeq{Cmaxlbd1},
we can choose $ \constt$ and $V^{-1}$ sufficiently small that, say,
\eq\lbeq{Cmax-low-a}
    \Exp\big(|\Cmax|\big)
    \geq
    \frac{\chi^2(p)}{3600} \Pro\Big(|\Cmax|\geq \frac{\chi^2(p)}{3600}\Big)
    \geq
    10^{-4}\chi^2(p).
\en
This gives the lower bound of \refeq{cmaxbd1}.
%
\qed

\section{Upper bound on the  supercritical susceptibility}
\label{sec-chi-sup}

In this section, we prove the bound \refeq{chiasysup}
of Theorem~\ref{main-thm-sup}~ii), by showing that
if $\constt$ and $\constoplow$ are sufficiently
small and $p=p_c+\cn^{-1}\epsilon$ with $0\leq\epsilon\leq 1$, then
\eq
\lbeq{chi-sup-up-final}
    \chi(p)\leq 81 V^{1/3}+(81\epsilon)^2V.
\en
The proof of \refeq{chi-sup-up-final} is based on the decomposition
\eq\lbeq{chisum}
\chi(p)
=\chi(p,\gamma)+\chi_\perp(p,\gamma)
\en
discussed in Section~\ref{sec-sucov}, where
\eq
\lbeq{chiMpg-rew}
    \chi(p,\gamma) = \Expg\big(|C(0)|I(0\not\leftrightarrow \Gcal)\big)
    =
    \sum_{k=0}^{V} k (1-\gamma)^k \Pbold_{p} (|C(0)|=k)
\en
and
\eq
\lbeq{chiperp-rew}
\chi_\perp(p,\gamma)=
\Expg\big(|C(0)|I(0\leftrightarrow \Gcal)\big)
=
    \sum_{k=0}^{V} k[1-(1-\gamma)^k] \Pbold_{p} (|C(0)|=k).
\en

For an upper bound on $\chi(p,\gamma)$, it follows from
Proposition~\ref{Mleq} and the lower  bound of \refeq{M-gamma-chi} that
\eq\lbeq{chimag-up1}
        \chi(p,\gamma)\leq
    \sqrt{\frac{12}\gamma}+\frac{13\epsilon} {\gamma}.
    \en
whenever $\constt$ and $\constoplow$ are sufficiently
small and $p=p_c+\cn^{-1}\epsilon\geq p_c$.  This gives a bound
$O(\epsilon^{-1})$, if we choose
 $\gamma$ proportional to $\epsilon^2$.
To obtain a bound of the form $O(\epsilon^2 V)$
for $\chi_\perp(p,\epsilon^2)$, we will make use of the random
variable
$Z_{\Gcal}=
\sum_{x\in\ver}I(x\leftrightarrow \Gcal)$.
As a first step, we prove the following two lemmas,
which give bounds on $\chi_\perp(p,\gamma)$ that are valid for all $p$
and $\gamma$.

\begin{lemma}
\label{lem:chiperpZ}
Let $0\leq p \leq 1$ and $0\leq\gamma\leq 1$.
Then
\eq\lbeq{chiperp-bd1}
\chi_\perp(p,\gamma)\leq
\frac 1V
\Expg(Z_{\Gcal}^2)
\leq V M^2(p,\gamma)+\chi_\perp(p,\gamma),
\en
\end{lemma}

\proof
Under the condition that $0\leftrightarrow\Gcal$,
$|C(0)|$ can be bounded by the number of vertices
that are connected to a green vertex, so that
\eq
|C(0)| I(0\leftrightarrow\Gcal)\leq
\sum_{x\in\ver}I(x\leftrightarrow\Gcal)I(0\leftrightarrow\Gcal).
\en
Combined with transitivity and the definition of
$Z_{\Gcal}$, this implies the lower bound in
\refeq{chiperp-bd1}.

To prove the upper bound, we decompose
the expectation of $Z_{\Gcal}^2$ as
\eqalign
\Expg(Z_{\Gcal}^2)
=&\sum_{x,y\in\ver}
\Expg[I(x\leftrightarrow\Gcal)I(y\leftrightarrow\Gcal)
I(x\not\leftrightarrow y)]
\nnb
&+\sum_{x,y\in\ver}
\Expg[I(x\leftrightarrow\Gcal)I(y\leftrightarrow\Gcal)
I(x\leftrightarrow y)].
\enalign
As an upper bound, the three events in the first term can
be replaced by $\{x \conn \Gcal\}\circ \{ y \conn \Gcal\}$.
We then use the BK inequality (with respect to the joint
bond/vertex measure) to bound the first term
by
\eqalign
\sum_{x,y\in\ver}
\Expg[I(x\leftrightarrow\Gcal)]
\Expg[I(y\leftrightarrow\Gcal)]
=V^2 M^2(p,\gamma).
\enalign
Since the second term can be rewritten as
\eq
\sum_{x,y\in\ver}
\Expg[I(x\leftrightarrow\Gcal)
I(x\leftrightarrow y)]
=\sum_{x\in\ver}
\Expg[|C(x)|I(x\leftrightarrow\Gcal)]
=
V\chi_\perp(p,\gamma),
\en
this proves the upper bound in \refeq{chiperp-bd1}.
\qed

\begin{lemma}\label{lem:chiperpbd}
Let $0\leq p \leq 1$ and $0\leq\gamma\leq 1$.
Then
\eq\lbeq{chiperp-bd3}
\frac \gamma{1-\gamma}
\Expg\big[|C(0)|^2I(0\not\leftrightarrow \Gcal)\big]
\leq
\chi_\perp(p,\gamma)
\leq
\gamma\chi^3(p).
\en
\end{lemma}

\proof
We first note that
\eq
\frac 1{1-\gamma}
\Expg\big[|C(0)|^2I(0\not\leftrightarrow \Gcal)\big]
=
\sum_{k=1}^V(1-\gamma)^{k-1}k^2 P_k(p)
=
\frac{\partial\chi_\perp(p,\gamma)}{\partial\gamma}
\en
is monotone decreasing in $\gamma$.  Using this fact and the
observation that $\chi_\perp(p,0)=0$, integration over $[0,\gamma]$ gives
\eq
\gamma\frac 1{1-\gamma}
\Expg\big[|C(0)|^2I(0\not\leftrightarrow \Gcal)\big]
\leq\chi_\perp(p,\gamma)
\leq
\gamma\bigg[\frac 1{1-\gamma}
\Expg\big[|C(0)|^2I(0\not\leftrightarrow \Gcal)\big]\bigg]_{\gamma=0}
\en
The right hand side is simply $\gamma \Exp[|C(0)|^2]$,
which is less than $\gamma \chi(p)^3$ by the tree-graph inequalities
\cite{AN84}.
\qed

In Appendix~\ref{sec-Zsup}, we use
Lemmas~\ref{lem:chiperpZ}--\ref{lem:chiperpbd} to
derive the differential inequality
\eq\lbeq{der-bd-final}
\frac{\partial}{\partial p}
\Expg\big[Z^2_\Gcal\big]
\leq
\frac {3\cn}{1-p}\frac{1-\gamma}\gamma
M(p,\gamma)\Expg\big[Z^2_\Gcal\big],
\en
for any $0\leq p\leq 1$ and $0\leq\gamma\leq 1$.
We use this to prove the following
lemma, which is the final
ingredient needed for the proof of \refeq{chi-sup-up-final}.

\begin{lemma}
\label{lem:Zsup}
If $\constt$ and $\constoplow$ are sufficiently
small,
$p=p_c+\cn^{-1}\epsilon\geq p_c$ and $0\leq\gamma\leq 1$, then
\eq\lbeq{Zsup-bd1}
\chi_\perp(p,\gamma)
\leq
13\gamma V
\exp\Big\{
\frac {3({1-\gamma})}{1-p}
\Big(\sqrt{12\frac{\epsilon^2}\gamma}+13\frac{\epsilon^2}\gamma\Big)\Big\}.
\en
\end{lemma}

\proof
We divide \refeq{der-bd-final} by the expectation on its right side
and integrate over the interval $[p_c,p]$.  Since
$M(p,\gamma)/(1-p)$ is monotone increasing in $p$, the right side (after
the above division) is bounded by its value at the upper limit $p$
of integration.  This leads to
\eq\lbeq{ZGcal-bd1}
\Expg\big[Z^2_\Gcal\big]
\leq
\Expcg\big[Z^2_\Gcal\big]
\exp\Big\{
\frac {3\epsilon}{1-p}\frac{1-\gamma}\gamma
M(p,\gamma)\Big\}.
\en
We will show that
\eq
\lbeq{Z2pc}
\Expcg(Z_{\Gcal}^2)
\leq 13\gamma V^2.
\en
With \refeq{chiperp-bd1},
\refeq{ZGcal-bd1} and Proposition~\ref{Mleq}, this gives
the desired estimate.
To prove \refeq{Z2pc}, we combine the bounds of Lemmas~\ref{lem:chiperpZ}
and \ref{lem:chiperpbd}
with Proposition~\ref{Mleq},
to get
\eq
\Expcg(Z_{\Gcal}^2)
\leq 12 \gamma V^2 +\gamma\chi^3(p_c)V
    = 12 \gamma V^2 +\gamma\lambda^3 V^2.
\en
It suffices to show that $\lambda \leq 1$.
If $V=1$, this bound is trivial, so let us assume
that  $V\geq 2$.  But in this case we may use
the bound \refeq{lambda-bd} to conclude that
$\lambda\leq 1$ whenever $\constt\leq 1/2$.
\qed

\medskip \noindent
{\em Proof  of Theorem~\ref{main-thm-sup}~ii).}
It follows from \refeq{chisum}, \refeq{chimag-up1} and \refeq{Zsup-bd1}
that
\eq
    \chi(p)\leq
    \frac 1\epsilon
    \bigg(\sqrt{12\frac{\epsilon^2}\gamma}+13\frac{\epsilon^2} {\gamma}\bigg)
    +
    13\gamma V
    \exp\Big\{
    \frac {3}{1-p}
    \Big(\sqrt{12\frac{\epsilon^2}\gamma}+13\frac{\epsilon^2}\gamma\Big)\Big\}.
\en
To estimate the factor $1-p$, we use the bound \refeq{pc-less-constt}
and the definition \refeq{opc-1-low} of $\constoplow$
to see that
\eq
1-p=1-p_c-\cn^{-1}\epsilon
\geq 1-\constt -\frac{\constt}{1-\constoplow}
\en
whenever $0\leq\epsilon\leq 1$.  Setting $\gamma=52\epsilon^2$,
and assuming that $\constt$ and $\constoplow$ are chosen small enough
to ensure that
$\frac 1{1-p}(\sqrt{3/13}+1/4)\leq 3/4$, we then get
\eq
\lbeq{821}
    \chi(p)\leq
    \frac 3{4\epsilon}
        +
    (26\epsilon)^2 V e^{9/4}
    \leq
    \frac 3{4\epsilon}
        +
    (81\epsilon)^2 V.
\en
Let $\epsilon' = \frac{1}{81}V^{-1/3}$.
We distinguish the cases $\epsilon < \epsilon^\prime$
and $\epsilon\geq \epsilon^\prime$.  In the first case,
we use monotonicity and \refeq{821} to obtain
\eq
    \chi(p)\leq
    \chi(p_c+\cn^{-1}\epsilon^\prime)
    \leq
    \frac 34 81 V^{1/3}
            +
    V^{1/3}
    \leq 81 V^{1/3}.
\en
If $\epsilon \geq \epsilon^\prime$,
we use $\epsilon^{-1} \leq 81V^{1/3}$ to obtain
\eq
    \chi(p)\leq
     \frac 34 81 V^{1/3}
        +
    (81\epsilon)^2 V
     \leq 81 V^{1/3}
     +
    (81\epsilon)^2 V
    .
\en
The combination of these two estimates gives \refeq{chi-sup-up-final}.
\qed

\section{Upper bound on the largest supercritical cluster}
\label{sec-supper-upb}

In this section, we prove the upper bound on the largest
supercritical cluster stated in Theorem~\ref{main-thm-pp}~ii).
For $p \leq p_c$, we used the variance bound Lemma~\ref{var-lemma}.
For $p\geq p_c$, we will use the following alternate bound on the variance
of $Z_{\geq s}$.  For its statement, we define
\eq\lbeq{chils}
\chi_{< s}(p)=\Exp\left[|C(0)|I[|C(0)|<s\right].
\en

\begin{lemma}\label{var-lemma2}
Let $s >0$ and $p\in (0,1)$.
Then
\eq
\lbeq{Z-var-bound1}
\Var\left[Z_{\geq s}\right]
\leq
    (1+ p\,\cn  s) V\chi_{<s}(p)
    \leq 4 s(1+ p\,\cn  s)V M(p,s^{-1}).
\en
\end{lemma}

\begin{proof}
We define the random
variable $Z_{<s}=V-Z_{\geq s}= \sum_{v \in \mathbb V}I[|C(v)|<s]$
and express the variance as
    \eq\lbeq{varyeq}
    \Var[Z_{\geq s}]=\Var[Z_{<s}] = \sum_{v,w,S,T}
    \big[ \Pro[C(v)=S, C(w)=T] -
    \Pro[C(v)=S]\Pro[C(w)=T] \big],
    \en
where the sum is over connected sets $S,T$ with $|S|<s$,
$|T|<s$ and $v\in S$, $w\in T$.
Let $\text{dist}(\cdot,\cdot)$ denoting the graph
distance on $\gr$.  If $\dist(S,T)>1$, the
above events are independent and so the difference is zero.
If $\dist(S,T)=0$, then $S\cap T\neq 0$, and the
first probability is zero unless $S=T$.  The corresponding
contribution from the first term is just
\eq\lbeq{term-1}
\sum_{S:|S|<s,\atop v,w\in S}\Pro[C(v)=S]
=V\chi_{<s}(p),
\en
implying that the contribution from both terms
can be bounded by \refeq{term-1}.

We are left with the contribution of the terms
$(v,w,S,T)$ with $\dist(S,T)=1$.  We need some
notation.  Given a connected set $S$, let
$E(S)$ be the set of all edges with both endpoint
in $S$, and $\partial S$ be the set of edges
with exactly one endpoint in $S$.  Finally,
let $A(S)$ be the event that $E(S)$ contains
a set $E$ of occupied edges such that the graph $(S,E)$
is connected.  With this notation, the event
$C(v)=S$ is just the intersection of the event
$A(S)$ with the event that all edges in
$\partial S$ are vacant.  Note also, that the
two events are independent, so that
$\Pbold_p(C(v)=S)$ is the product
$\Pbold_p[A(S)](1-p)^{|\partial S|}$.

If $\dist(S,T)= 1$, the events
$A(S)$, $A(T)$ and the event that the
edges in $\partial S\cup\partial T$ are
vacant are independent.  For these
$(S,T)$, the
difference in \refeq{varyeq} can thus
be rewritten as
 \eqalign
     \Pro[C(v)=S, & C(w)=T] -
    \Pro[C(v)=S]\Pro[C(w)=T]
    \nonumber
    \\
&= \Pro[A(S)]\Pro[A(T)]
\Bigl((1-p)^{|\partial S\cup\partial T|}
-(1-p)^{|\partial S|}(1-p)^{|\partial T|}\Bigr)
\nonumber
\\
\lbeq{term2}
&=\Pro[C(v)=S, C(w)=T]\Bigl(1
-(1-p)^{|\partial S\cap\partial T|}\Bigr).
    \enalign
 To continue, we use the inequality
$1- (1-p)^k\leq p k$ to obtain
\eq\lbeq{bound-term2}
1-(1-p)^{-|\partial S\cap\partial T|}
\leq p|\partial S\cap\partial T|
=p\sum_{x,y:\dist(x,y)=1}I[x\in S]I[y\in T].
\en
Combining \refeq{bound-term2} with the identity
\refeq{term2}, we now
bound the contribution to \refeq{varyeq} due to the terms
$(v,w,S,T)$ with $\dist(S,T)=1$ by
    \eqalign
    p\sum_{v,w,x,y\atop\text{dist}(x,y)=1}
    \Pro&\big(|C(v)|<s,x\in C(v),|C(w)|<s,y\in C(w), C(x)\neq C(y)\big)
    \nonumber
    \\
    &=
    p\sum_{x,y,v,w\atop\text{dist}(x,y)=1}
    \Pro\big(|C(x)|<s,  v \in C(x), |C(y)|<s, w\in C(y),C(x)\neq C(y)\big)
    \nonumber
    \\
    &=
    p\sum_{x,y\atop\text{dist}(x,y)=1}
    \Exp\big(|C(x)|\, I[|C(x)|<s]|C(y)|\, I[|C(y)|<s,x\not\conn y]\big)
    \nonumber
    \\
    &\leq
     ps\sum_{x,y\atop\text{dist}(x,y)=1}
    \Exp\big(|C(x)|\, I[|C(x)|<s]\big)
    \nonumber
    \\
    &    = p s\cn  V\chi_{<s}(p).
    \enalign
Combining this term with \refeq{term-1}, we obtain the first bound of
\refeq{Z-var-bound1}.

For the second bound of \refeq{Z-var-bound1}, it suffices to show
that
\eq
\lbeq{chilM}
    \chi_{<s}(p) \leq 4sM(p,s^{-1}).
\en
For $s \geq 2$,  we bound $\chi(p,s^{-1})$ in \refeq{chiMpg} from
below by restricting the sum over $k$ to $k<s$.  We then use
$(1-s^{-1})^k\geq(1-s^{-1})^s\geq 1/4$ to conclude that
$\chi_{<s}(p)\leq 4\chi(p,s^{-1})$.
Combined with \refeq{M-gamma-chi}, this gives
\refeq{chilM} for $s \geq 2$.  If $s \leq 1$, the left side of \refeq{chilM}
is zero and the bound is trivial.  Finally, if $1<s<2$ then the left side
of \refeq{chilM} is $1$, whereas it follows from the fact that
$M(p,\gamma) \geq \Pbold_{p,\gamma} (0 \in \Gcal) = \gamma$ that
the right side of \refeq{chilM} is at least $4$.
\end{proof}


\smallskip \noindent
{\em Proof of Theorem~\ref{main-thm-pp}~ii).}
Let $p=p_c+\cn^{-1}\epsilon$.  It suffices to prove that under the
hypotheses of Theorem~\ref{main-thm-pp} there are constants $b_{11}$, $b_{12}$
such that
\eq\lbeq{cmax-up-precise}
    \Pro\Big(|\Cmax|\geq
        [1+(\epsilon V^\eta)^{-1}]
        V \theta_\alpha(p)\Big)
    \leq
    \frac {b_{11}}{(\epsilon V^\eta)^{3-2\alpha}}
\en
if $\max\{b_{12} V^{-1/3},V^{-\eta}\}\leq\epsilon\leq 1$.
(The proof actually applies for $b_{12} V^{-1/3} \leq\epsilon\leq 1$,
but the result is not meaningful unless $\epsilon \geq V^{-\eta}$.)
To prove \refeq{cmax-up-precise},
we first note that $V\theta_\alpha(p)\geq N_\alpha$
if $b_{10}[\epsilon V^{1/3}]^{3-\alpha}\geq 1$, by \refeq{bound1}.
To satisfy $b_{10}[\epsilon V^{1/3}]^{3-\alpha}\geq 1$, we will take
$b_{12} \geq b_{10}^{-1/(3-\alpha)}$.
Thus we have
\eqalign
    \Pro\Big(|\Cmax|\geq
        [1+(\epsilon V^\eta)^{-1}]
        V \theta_\alpha(p)\Big)
    &= \Pro\Big(Z_{\geq N_\alpha} \geq |\Cmax|\geq
        [1+(\epsilon V^\eta)^{-1}]
        V \theta_\alpha(p)\Big)
    \nnb
    &\leq
    \Pro\Big(|Z_{\geq N_\alpha} -
        V \theta_\alpha(p)| \geq (\epsilon V^\eta)^{-1}
    V \theta_\alpha(p)\Big).
\enalign
It therefore suffices to show that
if $\max\{b_{12} V^{-1/3},V^{-\eta}\}\leq\epsilon\leq 1$ then
\eq
\lbeq{Znear}
    \Pro\Big(|Z_{\geq N_\alpha} -
        V \theta_\alpha(p)| \geq (\epsilon V^\eta)^{-1}
    V \theta_\alpha(p)\Big)
    \leq
    \frac {b_{11}}{(\epsilon V^\eta)^{3-2\alpha}}.
\en

By the variance estimate \refeq{Z-var-bound1} and the triangle condition,
${\rm Var} [Z_{\geq s}] \leq 12 s^2 V M(p,s^{-1})$.
Since $M(p,\gamma) \leq \sqrt{12 \gamma}+13 \epsilon$ by \refeq{Mleq1},
it follows from the fact that $N_\alpha^{-1} \leq \epsilon^2$ that
\eqalign
    {\rm Var} [Z_{\geq N_\alpha}]
    &\leq
    200 \epsilon N_\alpha^2 V
    = 200
    \frac{1}{\epsilon^{3-2\alpha}} V^{1+2\alpha/3} .
\enalign
Therefore,
by Chebyshev's inequality  and \refeq{bound1},
\eq
    \Pbold \big(|Z_{\geq N_\alpha} - V\theta_\alpha| \geq
    (\epsilon V^\eta)^{-1} V\theta_\alpha \big)
    \leq
    (\epsilon V^\eta)^{2}
    \frac{{\rm Var} [Z_{\geq N_\alpha}]}{V^2\theta_\alpha^2}
    \leq
    \frac{200}{b_{10}^2}
    \frac{1}{\epsilon^{3-2\alpha}V^{1-2\eta - 2\alpha/3}}.
\en
Since $\eta = \frac{3-2\alpha}{15-6\alpha}$, the important factor on
the right side is equal to $(\epsilon V^\eta)^{-(3-2\alpha)}$.
This gives \refeq{Znear} and completes the proof.
\qed

\renewcommand{\thesection}{\Alph{section}}
\setcounter{section}{0}

\section{Appendix: Derivation of differential inequalities}
\label{sec-app}

\subsection{Differential inequality for the susceptibility}
\label{sec-chidif}

In this section, we prove \refeq{dis}, which is restated here as
Proposition~\ref{prop-chi'}.  We follow the original
proof of Aizenman and Newman \cite{AN84}, with a minor extension
for the lower bound
to deal with an arbitrary transitive graph
$\gr$.  The proof also provides an instructive preliminary to the proof
of \refeq{rdi} in Appendix~\ref{sec-M-diff-ine}.

\begin{prop}
\label{prop-chi'}
For all $p\in (0,1)$,
    \eq
    \lbeq{chi'ul}
    [1-\bar\nabla_p] \cn\chi(p)^2 \leq
    \frac {d \chi(p)}{dp} \leq \cn\chi(p)^2.
    \en
\end{prop}

Recall that $E\circ F$
denotes the event that $E$ and $F$ occur disjointly.
Given a bond
configuration, we
    say that a bond is {\em pivotal}\/ for $x \conn y$ if
    $x \conn y$ in the possibly modified configuration in which
    the bond is made occupied, whereas $x$ is not connected to $y$
    in the possibly modified configuration in which
    the bond is made vacant.

\smallskip \noindent
{\em Proof of the upper bound in \refeq{chi'ul}.}
By Russo's formula (see \cite[Theorem~(2.25)]{Grim99}),
    \eq\lbeq{difftau}
    \frac{d}{dp} \tau_p(x,y)
    = \sum_{\{u,v\}\in \mathbb B}
    \Pro(\{u,v\} \text{ is pivotal for }x\conn y).
    \en
Therefore, by the BK inequality,
    \eq
    \frac{d}{dp} \tau_p(x,y)
    \leq \sum_{(u,v)} \Pro(\{x\conn u\} \circ \{v\conn y\})
    \leq \sum_{(u,v)} \tau_p(x,u)\tau_p(v,y),
    \en
where the sum over $(u,v)$ is a sum over {\em directed} bonds.
We then perform the sums over
$y,v,u$ (in that order) and use transitivity to obtain
the desired upper bound.
\qed

\smallskip
For the lower bound of \refeq{chi'ul}, we
will use the following definition
and lemmas.  In the first lemma, we use transitivity to give an
alternate representation for $\sum_{v : \{0,v\}\in \edg}\nabla_p(0,v)$.
This is related to an issue raised by Schonmann \cite{Scho01}
(see also \cite{Scho02}), who pointed out
that the use of differential inequalities
plus the triangle
condition to prove mean-field behavior on general {\em infinite} transitive
graphs can be accomplished under the additional assumption that the
graph is unimodular,
and that it is an open problem to determine whether the
assumption of unimodularity is essential.  Finite transitive
graphs are always unimodular, so the issue raised in
\cite{Scho01} is less relevant for our
purposes.  In any case, we will bypass the issue altogether
by applying the following lemma.  For its statement, we define
\eq
    T_1(z) =
    \sum_{(u,v)}\sum_{y \in \mathbb V} \tau_p(z,u)\tau_p(z,y)\tau_p(y,v).
\en
The equality \refeq{nouni} of Lemma~\ref{lem-nouni} will be applied only
in \refeq{dif2} and \refeq{X3ub}.

\begin{lemma}
\label{lem-nouni}
For each $u,z \in \mathbb V$,
\eq
\lbeq{nouni}
    T_1(z)
    =
    \sum_{v : \{u,v\}\in \edg} \nabla_p(u,v) \leq \cn \bar\nabla_p.
\en
\end{lemma}

\proof
The inequality follows from the definition of $\bar\nabla_p$ in
\refeq{nabbardef}.  To prove the equality, let
\eq
    T_2(u)
    =
    \sum_{v:\{u,v\}\in \edg} \nabla_p(u,v)
    = \sum_{v:\{u,v\}\in \edg}\sum_{y,z \in \mathbb V}
    \tau_p(u,z)\tau_p(z,y)\tau_p(y,v)
    .
\en
We first prove that $T_2(u)$ is independent of $u$; a similar proof applies
for $T_1(z)$.  By transitivity, there is a graph automorphism
$\varphi = \varphi_u$
such that $\varphi(u)=0$, where $0$ is a fixed vertex.
Since $\tau_p(x,y) = \tau_p(\varphi(x),\varphi(y))$,
\eqalign
    T_2(u)& =
    \sum_{v:\{u,v\}\in \edg}\sum_{y,z \in \mathbb V}
    \tau_p(0,\varphi(z))\tau_p(\varphi(z),\varphi(y))
    \tau_p(\varphi(y),\varphi(v)).
\enalign
Since $\varphi$ is an automorphism,
$\{u,v\}\in \edg$   if and only if
$\{\varphi(u),\varphi(v)\}=\{0,\varphi(v)\}\in \edg$.
Similarly, as $x$ runs over all vertices, so does $\varphi(x)$.  Relabelling
$\varphi(y)$ to $y$, $\varphi(v)$ to $v$, and $\varphi(z)$ to $z$, we thus get
\eq
    T_2(u)=
    \sum_{v:\{0,v\}\in \edg}\sum_{y,z \in \mathbb V}
     \tau_p(0,z)\tau_p(z,y)\tau_p(y,v)
      =T_2(0).
\en

Since $T_i(x)$ ($i=1,2$)
is independent of $x \in \mathbb V$,
it is equal to the average of its sum over $x\in \mathbb V$.
Since $\sum_{z \in \mathbb V}T_1(z)=\sum_{u \in \mathbb V}T_2(u)$,
this implies that $T_1(z) = T_2(u)$ for all $u,z \in \mathbb V$,
which proves the equality in \refeq{nouni}.
\qed

    \begin{defn}
    \label{def-inon}
    (a)
    Given a bond configuration, and $A \subset \mathbb V$, we
        say $x$ and $y$ are \emph{connected in $A$}, if there is an
        occupied path from $x$ to $y$ having all its endpoints in
        $A$, or if $x = y \in A$.
        We define a restricted two-point function by
            \eq
            \tau^{A} (x, y) = \prob{\text{$x$ and $y$ are
            connected in $\mathbb V\backslash A$}}
            .
            \en
        (b)
        Given a bond configuration, and $A\subset \mathbb V$, we
        say $x$ and $y$ are \emph{connected through $A$}, if
    $x\conn y$ and every
        occupied path connecting $x$ to $y$ has at least one bond
        with an endpoint in $A$.  This event is written as
        $x \ct{A} y$.
    \newline
        (c)
        Given a bond configuration, and a bond $b$, we define
        $\tilde{C}^{b}(x)$ to be the set of vertices connected to $x$,
        in the new configuration obtained by setting $b$ to be vacant.
        \newline
    (d) Given an event $E$, we define the event $\{E$ occurs on
    $\tilde{C}^{(u,v)}(x)\}$ to be the set of configurations such
    that $E$ occurs on the modified configuration in which every
    bond that does not have an endpoint in $\tilde{C}^{(u,v)}(x)$
    is made vacant. We say that $\{E$ occurs in $\mathbb V \backslash
    \tilde{C}^{(u,v)}(x)\}$ if $E$ occurs on the modified
    configuration in which every bond that does not have both
    endpoints in $\mathbb V \backslash \tilde{C}^{(u,v)}(x)$ is made
    vacant.
    \end{defn}

\begin{lemma}
    \label{lem-cut1}
    Fix $p \in [0,1]$.
    Given a bond $(u, v)$, a vertex $w$ and events $E, F$,
    \eqalign
    \lbeq{lemcut1}
        & \Ebold_p \left(
        I[E \text{{\rm~occurs on }} \tilde{C}^{(u, v)}(w)
        \AND
        (u,v) \text{\rm ~is occupied } \AND
        F \text{{\rm~occurs in }} \mathbb V\backslash
        \tilde{C}^{(u, v)}(w) ]
        \right)
        \nnb
        & \qquad
        =
        p
        \Ebold_p \left( I[E \text{{\rm~occurs on }}
        \tilde{C}^{(u, v)}(w)]
        \,
        \Ebold_p \left( I[F \text{{\rm~occurs in }} \mathbb V\backslash
        \tilde{C}^{(u, v)}(w) ]
        \right)
        \right) .
    \enalign
The identity \refeq{lemcut1} is also valid if the event $\{(u,v)$
is occupied$\}$ is removed from the left side and $p$ is
removed from the right side.
\end{lemma}

The above lemma is present in \cite{AN84} in an implicit form.
Its elementary proof can be found in \cite[Lemma~3.2]{BCHSS04b}.
In the nested expectation on the right side of \refeq{lemcut1},
the set $\tilde{C}^{(u,
v)}(w)$ is a random set with respect to the outer expectation,
but it is deterministic with respect to the inner expectation.  The
inner expectation on the right side effectively introduces a second
percolation model on a second lattice, which is coupled to the original
percolation model via the set $\tilde{C}^{(u, v)}(w)$.


\smallskip \noindent
{\em Proof of the lower bound in \refeq{chi'ul}.}
By definition,
    \eqalign
    \lbeq{rewriteonin}
    & \{(u,v) \text{ is pivotal for }0\conn x \}
    \\ \nonumber & \quad = \{0 \conn u \text{ occurs on }
    \tilde{C}^{(u, v)}(0)\}\cap
    \{v \conn x \text{ occurs in }\mathbb V \backslash \tilde{C}^{(u, v)}(0)\}.
    \enalign
By Lemma~\ref{lem-cut1} and \refeq{rewriteonin},
    \eqalign
    \Pro((u,v) \text{ is pivotal for }0\conn x)
    & =\Ebold \left( I[0 \conn u \text{ occurs on }\tilde{C}^{(u, v)}(0)]
    \,
    \tau^{\tilde{C}^{(u, v)}(0)}(v, x)
    \right)
    \nnb
    \lbeq{inondrop}
    & =
    \Ebold \left( I[0 \conn u ] \,
    \tau^{\tilde{C}^{(u, v)}(0)}(v, x)
    \right).
    \enalign
In the second equality of \refeq{inondrop},
we dropped the condition ``occurs on
$\tilde{C}^{(u, v)}(0)$,'' because of the fact that
$\tau^{\tilde{C}^{(u,v)}(0)}(v,x)=0$ on the event
$\{0 \conn u\}\setminus \{0 \conn u \text{ occurs on }
\tilde{C}^{(u, v)}(0)\}$.
The identity \refeq{inondrop} can be rewritten as
    \eq
    \Pro((u,v) \text{ is pivotal for }0\conn x)
    =\tau_p(0, u) \tau_p(v,x)
    -\Ebold \left( I[0 \conn u]
    \Pro(v \ctx{\tilde{C}^{(u, v)}(0)} x)\right).
    \en

By the BK inequality, for $A \subset \mathbb V$ we have
    \eqalign
    \Pro(v \ctx{A} x)&=\Pro
    \Big(\bigcup_{y\in A} \{v \conn y\}
    \circ \{y\conn x\}\Big) \nonumber \\
    &\leq
    \sum_{y\in \mathbb V} I[y\in A] \Pro(\{v \conn y\}
    \circ \{y\conn x\})\nonumber\\
    &\leq \sum_{y\in \mathbb V} I[y\in A]
    \tau_p(v,y)\tau_p(y,x).
  \lbeq{chidifbk}
    \enalign
Therefore, for $A = \tilde{C}^{(u, v)}(0) \subset C(0)$, we have
\eq
    \Pro(v \ctx{\tilde{C}^{(u, v)}(0)} x)
    \leq
    \sum_{y\in \mathbb V} I[y \in C(0)] \tau_p(v,y)\tau_p(y,x).
\en
Substitution yields
    \eq
    \Pro((u,v) \text{ is pivotal for }0\conn x)
    \geq
    \tau_p(0,u)\tau_p(v,x)
    -\sum_{y\in \mathbb V} \Pro(0 \conn u\conn y)\tau_p(v,y)\tau_p(y,x).
    \en
The tree-graph bound \cite{AN84} (which is an elementary consequence of
the BK inequality) implies that
    \eq
    \Pro(0 \conn u, 0 \conn y)\leq
    \sum_{z \in \mathbb V} \tau_p(0,z)\tau_p(z,y)\tau_p(z,u).
    \en
Therefore,
    \eqalign
    \Pro((u,v) \text{ is pivotal for }0\conn x)
    &\geq \tau_p(0,u)\tau_p(v,x)\\
    &\quad -
    \sum_{y,z \in \mathbb V}
    \tau_p(0,z)\tau_p(z,y)\tau_p(z,u)\tau_p(y,v)\tau_p(y,x).
    \nonumber
    \enalign
Recalling \refeq{difftau}, and performing the sums over $u,v,x$, leads to
    \eqalign
    \frac{d\chi(p)}{dp} &\geq  \cn\chi(p)^2-
    \chi(p)\sum_{z\in \mathbb V} \tau_p(0,z)\sum_{(u,v)}\sum_{y \in \mathbb V}
    \tau_p(z,y)\tau_p(z,u)\tau_p(y,v)\nnb
    &=
    \cn\chi(p)^2- \chi(p)\sum_{z\in \mathbb V}\tau_p(0,z) T_1(z)
    \nnb &
    \geq \cn \chi(p)^2 [1-\bar\nabla_p],
    \lbeq{dif2}
    \enalign
by Lemma~\ref{lem-nouni}.
\qed

\subsection{Differential inequality for the magnetization}
\label{sec-M-diff-ine}

In this section, we prove the differential inequality \refeq{rdi}.
Our method of proof for \refeq{rdi} is related to, but simpler than, the
method used in \cite{BA91} to prove an analogous statement for
percolation on $\Z^n$. See \cite[Section~3]{HS00a} for related results
for $\Z^n$ which are stronger but more difficult to prove.
We restate \refeq{rdi} as \refeq{rdiapp} in the following lemma.
Note that, by \refeq{chiMpg}, the factor
$(1-\gamma)  \partial M/\partial\gamma$
on the right side of \refeq{rdi} can be replaced by $\chi(p,\gamma)$.

    \begin{lemma}
    \label{ineq-app}
    If $0< p < 1$ and $0<\gamma <1$, then
       \eq
        \lbeq{rdiapp}
         M(p,\gamma) \geq
        \left[{\cn \choose 2}p^2(1-p)^{\cn -2}(1-\nabla_p^{\rm max})^3
        -p-\nabla_p^{\rm max} \right]p\cn(1-\gamma)
    M^2(p,\gamma) \frac{\partial M(p,\gamma)}{\partial\gamma},
        \en
    where $\nabla_p^{\rm max} = \sup_{x,y\in \mathbb V}\nabla_p (x,y)$.
    \end{lemma}

\begin{proof}
Recall the use of the ``green'' set $\Gcal$ discussed in
Section~\ref{sec-sucov}.
Let $\{ v \dbc \Gcal \}$ denote the event that there exist
$x,y \in \Gcal$, with $x \neq y$, such that there are disjoint
connections $v \conn x$ and $v \conn y$.
Let $F_{(u,v)}$ denote the event that the bond $(u,v)$
is occupied and pivotal for the connection from $0$ to $\Gcal$,
with $\{ v \dbc \Gcal \}$.   Let $F = \cup_{(u,v)}F_{(u,v)}$, and
note that the union is disjoint.
Since $0 \conn \Gcal$ when $F$ occurs,
$M = \Pbold(0 \conn \Gcal) \geq \Pbold(F)$,
and it suffices to prove that $\Pbold(F)$ is bounded below by
the right side of \refeq{rdiapp}.

For $x,y \in \mathbb V$, we define a ``green-free''
analogue of the two-point function by
\eq
\lbeq{taugdef}
    \tau_{p,\gamma}(x,y)
    = \Pbold_{p, \gamma}(x \conn y, x \nc \Gcal),
\en
so that
\eq
\lbeq{chitaug}
    \chi(p,\gamma) = \sum_{x \in \mathbb V} \tau_{p,\gamma}(0,x)
\en
and $\chi(p,0) = \chi(p)$.
Given a subset $A \subset \mathbb V$, we define
\eqalign
    \tau_{p,\gamma}^A(x,y)
    &= \Pbold_{p,\gamma} ( (x\conn y, x \nc \Gcal) \mbox{ in }
    \mathbb V\backslash A).
\enalign
When $\gamma \neq 0$, we extend the definition of
``occurs in'' and ``occurs on'' in Definition~\ref{def-inon} as in
\cite[Definition~2.2]{HS00a}.  In particular, we now say that
$E$ occurs on $A$ if, given a configuration, $E$ occurs on the new
configuration obtained by setting all bonds not touching $A$ to be vacant
and all vertices not in $A$ to be not green.
By definition of $F_{(u,v)}$, it can be seen by conditioning on
the set $\tilde{C}^{(u,v)}(v)$ that
\eqalign
\lbeq{iii}
    \Pbold_{p,\gamma}(F)  & =
    p \sum_{(u,v)} \Pbold_{p,\gamma}
    \left[
    (0 \conn u \AND 0 \nc \Gcal) \mbox{ in }
    \mathbb V\backslash \tilde{C}^{(u,v)}(v),
    \;\; v\dbc \Gcal \mbox{ on } \tilde{C}^{(u,v)}(v)
    \right].
\enalign
It then follows from \cite[Lemma~2.4]{HS00a}
(a straightforward extension of Lemma~\ref{lem-cut1} to allow for the
presence of a magnetic field) that
\eqalign
    \Pbold_{p,\gamma}(F)  & =
    p \sum_{(u,v)} \Ebold_{p,\gamma}
    \left[
    \tau_{p,\gamma}^{\tilde{C}^{(u,v)}(v)}(0,u)
    I[v\dbc \Gcal \mbox{ on } \tilde{C}^{(u,v)}(v)]
    \right].
\enalign

We use the identities
\eq
    \tau_{p,\gamma}^{\tilde{C}^{(u,v)}(v)}(0,u)
    =
    \tau_{p,\gamma}(0,u) - \left(\tau_{p,\gamma}(0,u) -
    \tau_{p,\gamma}^{\tilde{C}^{(u,v)}(v)}(0,u)\right)
\en
and
\eq
    I[v\dbc \Gcal \mbox{ on } \tilde{C}^{(u,v)}(v)]
    =
    I[v\dbc \Gcal ]- \left( I[v\dbc \Gcal ] -
    I[v\dbc \Gcal \mbox{ on } \tilde{C}^{(u,v)}(v)] \right).
\en
Recalling \refeq{chitaug}, it follows that
\eqalign
\lbeq{iiib}
    P_{p,\gamma}(F) = \;\;
    & p\cn \chi(p,\gamma) \Pbold_{p,\gamma}(0 \dbc \Gcal)
    -
    p \sum_{(u,v)} \tau_{p,\gamma}(0,u) \Ebold_{p,\gamma}
    \left[
    I[v\dbc \Gcal ]
    - I[v\dbc \Gcal \mbox{ on } \tilde{C}^{(u,v)}(v)]
    \right]
    \nonumber \\ &
    -
    p \sum_{(u,v)} \Ebold_{p,\gamma}
    \left[
    \left(\tau_{p,\gamma}(0,u) -
    \tau_{p,\gamma}^{\tilde{C}^{\{u,v\}}(v)}(0,u)\right)
    I[v\dbc \Gcal \mbox{ on } \tilde{C}^{(u,v)}(v)]
    \right].
\enalign
We write \refeq{iiib} as $X_1-X_2-X_3$, bound $X_1$ from below,
and bound $X_2$ and $X_3$ from above.

\medskip \noindent
{\em Lower bound on $X_1$.} We will prove that
\eq
\lbeq{0dbcGbd}
    \Pbold_{p,\gamma}( 0 \dbc \Gcal)
    \geq
    {\cn \choose 2} p^2(1-p)^{\cn -2} M^2(p,\gamma)
    (1-\nabla_p^{\rm max})^3,
\en
which implies that
\eq
\lbeq{X1lb}
    X_1 \geq p\cn\chi(p,\gamma)
    {\cn \choose 2} p^2(1-p)^{\cn -2} M^2(p,\gamma) (1-\nabla_p^{\rm max})^3.
\en
 To begin, we note that the
event $\{ 0 \dbc \Gcal\}$ contains the event $\cup_{e,f}E_{e,f}$,
where the union is over unordered pairs of neighbors $e,f$ of the
origin, the union is disjoint, and the event $E_{e,f}$ is defined
as follows.  Let $E_{e,f}$ be the event that the bonds $(0,e)$ and
$(0,f)$ are occupied, all other bonds incident on $0$ are vacant,
and that in the reduced graph $\gr^-=(\mathbb V^-, \mathbb B^-)$
obtained by deleting the
origin and each of the $\cn$ bonds incident on $0$ from $\gr$ the
following three events occur:  $e \conn \Gcal$, $f \conn
\Gcal$, and $C(e) \cap C(f) = \varnothing$.  Let
$\Pbold_{p,\gamma}^-$ denote the joint bond/vertex measure on
$\gr^-$.  Then
    \eqalign
    \Pbold_{p,\gamma}(0 \dbc \Gcal)
    &
    \geq \Pbold_{p,\gamma} \big( \cup_{\{e,f\}} E_{e,f} \big)
    =
    \sum_{\{e,f\}}\Pbold_{p,\gamma}(E_{e,f})
    \nnb
    \lbeq{pdbG1}
    &
    = p^2(1-p)^{\cn -2} \sum_{\{e,f\}}
    \Pbold_{p,\gamma}^-( e \conn \Gcal, \;\; f \conn \Gcal, \;\;
    C(e) \cap C(f) = \varnothing).
\enalign

Let $W$ denote the event whose probability appears on the right
side of \refeq{pdbG1}.  Conditioning on the set $C(e)=A \subset
\mathbb V^-$, we see that
    \eq
    \Pbold_{p,\gamma}^-( W)
    =
    \sum_{A : A \ni e} \Pbold_{p,\gamma}^-( C(e)=A, \;\;
    e \conn \Gcal, \;\; f \conn \Gcal, \;\;
    C(e) \cap C(f) = \varnothing).
\en
This can be rewritten as
\eqalign
\lbeq{PW}
    \Pbold_{p,\gamma}^-( W)
    & =
    \sum_{A : A \ni e} \Pbold_{p,\gamma}^-( (C(e)=A, \;\;
    e \conn \Gcal) \mbox{  on   } A, \;\; f \conn \Gcal
    \mbox{  in  } \mathbb V^-\setminus A)
    \nnb
    & =
    \sum_{A : A \ni e} \Pbold_{p,\gamma}^-( C(e)=A, \;\;
    e \conn \Gcal) \; \Pbold_{p,\gamma}^-(f \conn \Gcal
    \mbox{  in  } \mathbb V^-\setminus A).
\enalign
Let $M^-(x) = \Pbold_{p,\gamma}^-(x \conn \Gcal)$, for $x \in \mathbb V^-$.
Then, by the BK inequality and the fact that the two-point function on
$\gr^-$ is bounded above by the two-point function on $\gr$,
\eq
\lbeq{PWlb}
    \Pbold_{p,\gamma}^-(f \conn \Gcal \mbox{  in  } \mathbb V^-\setminus A)
    = M^-(e) - \Pbold_{p,\gamma}^-(f \ct{A} \Gcal)
    \geq
    M^-(e) - \sum_{y \in A}\tau_{p,0}(f,y)M^-(y).
\en
By definition and the BK inequality,
\eq
\lbeq{M*}
    M^-(x) = M(p,\gamma) - \Pbold_{p,\gamma}(x \ctx{\{0\}} \Gcal)
    \geq
    M(p,\gamma)( 1 - \tau_{p,0}(0,x)) \geq M(p,\gamma)(1-\nabla_p^{\rm max}).
\en
In the above, we also used $\tau_{p,0}(0,x) \leq \nabla_p(0,x)$,
which follows from \refeq{tria-def} (with $u=v=y=0$).

It follows from \refeq{PW}--\refeq{M*} that
\eqalign
    \Pbold_{p,\gamma}^-( W)
    & \geq
    M(p,\gamma)(1-\nabla_p^{\rm max}) \sum_{A : A \ni e}
    \Pbold_{p,\gamma}^-( C(e)=A, \;\;
    e \conn \Gcal)
    \big[ 1 - \sum_{y \in A}\tau_{p,0}(f,y) \big]
    \nnb
\lbeq{M**}
    & =
    M(p,\gamma)(1-\nabla_p^{\rm max}) \big[ M^-(e)
    - \sum_{y \in \mathbb V^-} \tau_{p,0}(f,y)
    \Pbold_{p,\gamma}^-( e\conn y, \; e \conn \Gcal)\big].
\enalign
By the BK inequality,
\eq
    \Pbold_{p,\gamma}^-( e\conn y, \; e \conn \Gcal)
    \leq
    \sum_{w \in \mathbb V^-} \tau_{p,0}(e,w) \tau_{p,0}(w,y) M^-(w),
\en
and hence, by \refeq{M*}--\refeq{M**},
\eqalign
    \Pbold_{p,\gamma}^-( W)
    & \geq
    M(p,\gamma)(1-\nabla_p^{\rm max}) \big[ M^-(e)
    - \sum_{y,w \in \mathbb V^-} \tau_{p,0}(f,y)
    \tau_{p,0}(e,w) \tau_{p,0}(w,y) M^-(w) \big]
    \nnb
    & \geq
    M^2(p,\gamma)(1-\nabla_p^{\rm max})^3.
\enalign
This completes the proof of \refeq{0dbcGbd}, and hence of
\refeq{X1lb}.

\medskip \noindent
{\em Upper bound on $X_2$.}
This is the easiest term.  By definition,
\eq
    X_2
    =
    p \sum_{(u,v)} \tau_{p,\gamma}(0,u) \Ebold_{p,\gamma}
    \left[
    I[v\dbc \Gcal ]
    - I[v\dbc \Gcal \mbox{ on } \tilde{C}^{(u,v)}(v)]
    \right].
\en
For the difference of indicators to be nonzero, the double connection
from $v$ to $\Gcal$ must be realized via the bond $\{u,v\}$,
which therefore must be occupied.
The difference of indicators is therefore bounded above by the
indicator that the events $\{v \conn \Gcal\}$, $\{u \conn \Gcal\}$
and $\{\{u,v\} \text{ occupied}\}$ occur disjointly.
Thus we have
\eq
     \Ebold_{p,\gamma}
    \left[
    I[v\dbc \Gcal ]
    - I[v\dbc \Gcal \mbox{ on } \tilde{C}^{(u,v)}(v)]
    \right]
    \leq
    pM^2(p,\gamma),
\en
and hence
\eq
\lbeq{X2ub}
    X_2
    \leq
    p^2 \cn M^2(p,\gamma)\chi(p,\gamma) 
    .
\en

\medskip \noindent
{\em Upper bound on $X_3$.}
By definition,
\eq
\lbeq{X3-1}
    X_3
    =
    p \sum_{(u,v)} \Ebold_{p,\gamma}
    \left[
    \left(\tau_{p,\gamma}(0,u) -
    \tau_{p,\gamma}^{\tilde{C}^{(u,v)}(v)}(0,u)\right)
    I[v\dbc \Gcal \mbox{ on } \tilde{C}^{(u,v)}(v)]
    \right].
\en
The difference of two-point functions is the expectation of
\eqalign
    & I [ 0\conn u, 0 \nc \Gcal]
    -
    I [ 0\conn u \mbox{ in } \mathbb V\backslash \tilde{C}^{(u,v)}(v),
    0 \nc \Gcal]
    \nonumber \\
    & +
    I [ 0\conn u \mbox{ in } \mathbb V\backslash \tilde{C}^{(u,v)}(v),
    0 \nc \Gcal]
    - I[ (0\conn u, 0 \nc \Gcal))
    \mbox{ in } \mathbb V\backslash \tilde{C}^{(u,v)}(v) ]
    \nonumber \\
    &
    \hspace{7mm}
    \leq I[ 0 \ctx{\tilde{C}^{(u,v)}(v)} u, 0 \nc \Gcal],
\enalign
since the second line is non-positive
and the first line equals the third line.
Since the indicator in \refeq{X3-1} is bounded above by
$I[v\dbc \Gcal ]$, it follows that
\eq
\lbeq{X3z}
    X_3 \leq
    p \sum_{(u,v)} \Ebold_{p,\gamma}
    \left[
    \Pbold_{p,\gamma}
    [ 0 \ctx{\tilde{C}^{(u,v)}(v)} u, 0 \nc \Gcal]
    \;
    I[v\dbc \Gcal ]
    \right].
\en

By \cite[Lemma~4.3]{HS00a} (which is proved by conditioning on $\Gcal$),
\eq
\lbeq{oneh}
    \Pbold_{p,\gamma}
    [ 0 \ct{A} u, 0 \nc \Gcal]
    \leq
    \sum_{y\in \mathbb V} \tau_{p,\gamma} (0,y)\tau_{p,0}(y,u) I[y \in A].
\en
The important point in \refeq{oneh} is that the condition $0 \nc \Gcal$
on the left side is retained in the factor $\tau_{p,\gamma}(0,y)$ on the
right side (but not in $\tau_{p,0}(y,u)$).  With \refeq{X3z}, this gives
\eq
    X_3 \leq
    p \sum_{(u,v)} \sum_{y\in \mathbb V} \tau_{p,\gamma}(0,y)\tau_{p,0}(y,u)
    \Ebold_{p,\gamma}
    \left[
    I[v\dbc \Gcal ]I[y \in \tilde{C}^{(u,v)}(v)]
    \right].
\en
Since
\eq
    I[v\dbc \Gcal ]I[y \in \tilde{C}^{(u,v)}(v)]
    \leq
    I[\{v \conn w \conn y, v \conn \Gcal\} \circ \{v \conn \Gcal\}],
\en
a further application of BK gives
\eqalign
    X_3
    \lbeq{X3ub}
    & \leq
    p \sum_{y\in \mathbb V} \tau_{p,\gamma}(0,y) \sum_{(u,v)} \tau_{p,0}(y,u)
    \sum_{w\in \mathbb V} \tau_{p,0}(v,w) \tau_{p,0}(y,w) M^2(p,\gamma)
    \nnb & =
    p M^2(p,\gamma) \chi(p,\gamma) T_1(0) \leq
    pM^2(p,\gamma) \chi(p,\gamma) \cn \bar\nabla_p
    ,
\enalign
using Lemma~\ref{lem-nouni} in the last step.

The combination of \refeq{X1lb}, \refeq{X2ub} and
\refeq{X3ub} completes the proof of \refeq{rdiapp}.
\end{proof}

\subsection{The differential inequality \refeq{der-bd-final}}
\label{sec-Zsup}

Let $0\leq p\leq 1$, $0\leq\gamma\leq 1$,
and let $Z_\Gcal$ denote the number of vertices that are connected
to a green vertex.
The differential inequality \refeq{der-bd-final} states
that
\eq\lbeq{der-bd-final-app}
\frac{\partial}{\partial p}
\Expg\big[Z^2_\Gcal\big]
\leq
\frac {3\cn}{1-p}\frac{1-\gamma}\gamma
M(p,\gamma)\Expg\big[Z^2_\Gcal\big].
\en

\smallskip \noindent
{\em Proof of \refeq{der-bd-final-app}.}
Let $A_{x,y}$ be the event that $x\leftrightarrow\Gcal$ and
$y\leftrightarrow\Gcal$.  Then
\eq
\Expg\big[Z^2_\Gcal\big]=\sum_{x,y\in\ver}\Prog(A_{x,y}),
\en
and hence, by Russo's formula,
\eqalign
\frac{\partial}{\partial p}
\Expg\big[Z^2_\Gcal\big]
&=
\sum_{x,y\in\ver}
\sum_{\{u,v\}\in\edg}
\Prog(\{u,v\}\text{ is pivotal for }A_{x,y})
\\
&=\frac 1{1-p} \sum_{x,y\in\ver} \sum_{\{u,v\}\in\edg}
\Prog(
\{u,v\}\text{ is vacant and pivotal for }A_{x,y}).
\enalign
If $\{u,v\}$ is vacant and pivotal for
$A_{x,y}$, then exactly one of the two endpoints of the edge
$\{u,v\}$ is connected to a green vertex.  Moreover, if one of the two
endpoints of $\{u,v\}$ is connected to a green vertex, and
the other is not, then the edge $\{u,v\}$ is automatically vacant.
As a consequence,
\eq\lbeq{der-bd-1}
\frac{\partial}{\partial p}
\Expg\big[Z^2_\Gcal\big]
=\frac 1{1-p} \sum_{x,y\in\ver} \sum_{(u,v)}
\Prog(\{\{u,v\}\text{ is pivotal for }A_{x,y}\}
\cap
\{u\leftrightarrow\Gcal\}
\cap
\{v\not\leftrightarrow\Gcal\}
),
\en
where the sum over $(u,v)$ is a sum over {\em directed} edges.
To analyze the probability in \refeq{der-bd-1}, we distinguish
two cases: either exactly one of the two vertices $x$ and $y$
is connected to a green vertex, or neither of them is connected
to a green vertex.  It is not possible that both are connected to $\Gcal$,
because we are in a situation where $\{u,v\}$ is vacant,
and it cannot then also be pivotal for $A_{x,y}$.

Let us first estimate the contribution due to the event
that neither $x$ nor $y$ is connected to a green vertex.
A moment's reflection shows that this contribution can be rewritten
as
\eq
\lbeq{notxy}
\Prog(
\{u\leftrightarrow\Gcal\}
\cap
\{x\leftrightarrow y\leftrightarrow v\not\leftrightarrow\Gcal\}
).
\en
We will estimate \refeq{notxy} by applying the BK inequality, as generalized
by van den Berg and Fiebig \cite{BF87}
to cover intersections of increasing
and decreasing events, to the joint distribution
$\Pbold_{p,\gamma}$
(alternatively,
the decoupling
inequalities of \cite{BC96} could be applied).
With respect to $\Pbold_{p,\gamma}$, the event
$\{u\leftrightarrow\Gcal\}$ is increasing, whereas the event
$\{x\leftrightarrow y\leftrightarrow v\not\leftrightarrow\Gcal\}$ is
the intersection of an increasing and a decreasing event.  In addition,
these events must occur disjointly.
Therefore, by the BK inequality,
\refeq{notxy} is bounded by
\eq\lbeq{pivotbd1}
\Prog(u\leftrightarrow\Gcal)
\Prog(x\leftrightarrow y\leftrightarrow v\not\leftrightarrow\Gcal).
\en

Consider now the contribution from the event that
$x$ is connected to a green vertex, while $y$ is not.
This contribution can be rewritten as
\eq
\Prog(
\{u\leftrightarrow\Gcal\}
\cap
\{x\leftrightarrow\Gcal\}
\cap
\{y\leftrightarrow v\not\leftrightarrow\Gcal\}
),
\en
which we bound by
\eq\lbeq{pivotbd2}
\Prog(
\{u\leftrightarrow\Gcal\}
\cap
\{x\leftrightarrow\Gcal\})
\Prog(y\leftrightarrow v\not\leftrightarrow\Gcal).
\en
Interchanging the role of $x$ and $y$, we obtain
a similar bound on the contribution of the term
with $y\leftrightarrow\Gcal$ and
$x\not\leftrightarrow\Gcal$.

Inserting these three bounds into \refeq{der-bd-1}, and recalling
\refeq{chitaug}, we get
\eqalign
\frac{\partial}{\partial p}
\Expg\big[Z^2_\Gcal\big]
&\leq\frac 1{1-p} \sum_{x,y\in\ver} \sum_{(u,v)}
\Prog(u\leftrightarrow\Gcal)
\Prog(x\leftrightarrow y\leftrightarrow v\not\leftrightarrow\Gcal)
\notag
\\
&\quad +\frac 2{1-p} \sum_{x,y\in\ver} \sum_{(u,v)}
\Prog(
\{u\leftrightarrow\Gcal\}
\cap
\{x\leftrightarrow\Gcal\})
\Prog(y\leftrightarrow v\not\leftrightarrow\Gcal)
\notag
\\
& =
\frac {V\cn}{1-p}
M(p,\gamma)
\Expg[|C(0)|^2I(0\not\leftrightarrow\Gcal)]
+\frac {2\cn}{1-p}\chi(p,\gamma)
\Expg\big[Z_\Gcal^2\big].
\lbeq{der-bd-2}
\enalign
To complete the proof of \refeq{der-bd-final-app}, we estimate
\refeq{der-bd-2} by using Lemmas~\ref{lem:chiperpbd} and \ref{lem:chiperpZ}
for the first term, and using
the lower bound of \refeq{M-gamma-chi} for the second term.
\qed

\section*{Acknowledgments}
This work began during a conversation at afternoon tea,
while RvdH, GS and JS were visiting Microsoft Research.
The authors thank Benny Sudakov for bringing the question
of the critical point of the $n$-cube to our attention, and
for telling us about the recent papers \cite{ABS02} and \cite{FKM02}.
The work of GS was supported in part by NSERC of Canada.
The work of RvdH was carried out in part at the University
of British Columbia and in part at Delft University of Technology.

\bibliographystyle{plain}

\end{document}

%% file: def99c.tex
\def\UseSection{
        \numberwithin{equation}{section}
    \theoremstyle{plain}
        \newtheorem{theorem}    {Theorem}[section]
        \DefineTheorems 
}

\def\DefineTheorems{
    
    \newtheorem{lemma}      [theorem] {Lemma}
    
    \newtheorem{prop}       [theorem] {Proposition}
    
    \newtheorem{cor}        [theorem] {Corollary}

    \theoremstyle{definition}
    \newtheorem{defn}       [theorem] {Definition}

    \theoremstyle{definition}

}

\newcommand{\bt}   {\begin{theorem}}
\newcommand{\et}   {\end  {theorem}}
\newcommand{\bl}   {\begin{lemma}}
\newcommand{\el}   {\end  {lemma}}
\newcommand{\bp}   {\begin{prop}}
\newcommand{\ep}   {\end  {prop}}
\newcommand{\bc}   {\begin{cor}}
\newcommand{\ec}   {\end  {cor}}
\newcommand{\bd}   {\begin{defn}}
\newcommand{\ed}   {\end  {defn}}

\newcommand{\ba}   {\begin{array}}
\newcommand{\ea}   {\end  {array}}
\newcommand{\be}   {\begin{enumerate}}
\newcommand{\ee}   {\end  {enumerate}}
\newcommand{\bi}   {\begin{itemize}}
\newcommand{\ei}   {\end  {itemize}}

\def\eq#1\en{\begin{equation}#1\end{equation}}  
\def\eqsplit#1\ensplit{
    \begin{equation}\begin{split}#1\end{split}\end{equation}
    }
\def\eqalign#1\enalign{
    \begin{align}#1\end{align}
    }
\def\eqmul#1\enmul{
    \begin{multline}#1\end{multline}
    }
\newcommand{\eqarrstar} {\begin{eqnarray*}} 
\newcommand{\enarrstar} {\end{eqnarray*}} 
\newcommand{\eqarray}   {\begin{eqnarray}} 
\newcommand{\enarray}   {\end{eqnarray}} 
\newcommand{\nnb}   {\nonumber \\} 

\newcommand{\lbeq}[1]  {\label{e:#1}}
\newcommand{\refeq}[1] {\eqref{e:#1}}    

\makeatletter
\newcommand{\labelcounter}[2]{{%
    \stepcounter{#1}
    \protected@write\@auxout{}%
    {\string\newlabel{#2}{{\csname the#1\endcsname}{\thepage}}}%
    {\ref{#2}}
    }}
\makeatother
\newcommand{\Ebold} {{\mathbb E}}

\newcommand{\Pbold} {{\mathbb P}}

\newcommand{\Zbold} {{\mathbb Z}}

\newcommand{\Ccal}   {\mathcal{C}}

\newcommand{\Gcal}   {\mathcal{G}}

\newcommand{\spose}[1] {{\hbox to 0pt{#1\hss}} }
\newcommand{\ltapprox} {\mathrel{\spose{\lower 3pt\hbox{$\mathchar"218$}}
 \raise 2.0pt\hbox{$\mathchar"13C$}}}
\newcommand{\gtapprox} {\mathrel{\spose{\lower 3pt\hbox{$\mathchar"218$}}
 \raise 2.0pt\hbox{$\mathchar"13E$}}}

\newcommand{\dist} {{  \rm dist }}

